\documentclass[, 10 pt, conference]{ieeeconf}  

\IEEEoverridecommandlockouts                              

\usepackage{graphics} 
\usepackage{epsfig} 
\usepackage{mathptmx} 
\usepackage{times} 
\usepackage{amsmath} 
\usepackage{amssymb}  

\usepackage{url}
\usepackage[edges]{forest}
\let\labelindent\relax
\usepackage{enumitem}
\usepackage{caption}
\usepackage{subcaption}
\usepackage {xcolor}
\usepackage{setspace} 
\usepackage{algorithm}
\usepackage{algpseudocode}
\makeatletter
\newlength{\phaserulewidth}
\newcommand{\setphaserulewidth}{\setlength{\phaserulewidth}}

\setphaserulewidth{.5pt}

\title{\LARGE \bf
Learning to cooperatively estimate road surface friction}

\author{Jens-Patrick Langstrand and Maben Rabi$^{1}$
\thanks{$^{1}$Jens-Patrick Langstrand and Maben Rabi is with the Faculty of Computer Science, Engineering and Economics,
        {\O}stfold University College, Braveien 4, Halden, Norway
        {\tt\small jens.p.langstrand@hiof.no, maben.rabi@hiof.no}}%
}

\begin{document}

\maketitle
\thispagestyle{empty}
\pagestyle{empty}

\begin{abstract}
Under wintry conditions at narrow, curved mountainous roads, there is a higher risk of accidents. If we could track and report the current coefficients of friction between tyres and the different road surfaces, then we can reduce this risk. In specific, estimating  and reporting the friction experienced by vehicles that recently passed a given road section can help to warn following vehicles. %
To keep costs down, a potential module for friction estimation must be based on the standard  sensors already installed, such as the IMU, and sensors for the steering angle and wheel speeds. 
Existing algorithms do not satisfactorily estimate the friction, while  using only these sensors. 
For this, we propose a distributed system consisting of: (i)~processing of measurements from existing vehicular sensors, to implement a virtual sensor that captures the effect of low friction on 
the vehicle, (ii)~transmitting short  kinematic summaries from vehicles to a road side unit~(RSU), using V2X communication, and (iii)~estimating the 
friction coefficients,  
by running a machine learning regressor at the RSU, on summaries from individual vehicles, and then combining several such estimates.
      
In designing and implementing our system over  a road network, we face two key questions: (i)~should  each individual road section have a {\textit{local}} friction coefficient regressor, 
or can we use a {\textit{global}} regressor that covers all the possible road sections? and (ii)~how accurate are the resulting regressor estimates?
We test the performance of design variations of our solution, using simulations on the commercial package~Dyna4. We consider a single vehicle type with varying levels of tyre wear, and a range of road friction coefficients. We find that: (a)~only a marginal loss of accuracy is incurred in using  a global regressor as compared to local regressors,  
(b)~the consensus estimate at the RSU has a worst case error of about ten percent, if the combination is based on at least fifty recently passed vehicles, and
(c)~our regressors have root mean square~(RMS) errors that are less than five percent. The RMS error rate of our system is half as that of a commercial friction estimation service~\cite{andreasson2017luleaBachelorThesisNiraDynamicsMiscVersion}. 

But when tested with data from extreme driving manoeuvres that were unseen in the training data,  our regressor performs an order of magnitude worse than on data from normal driving runs on curved road sections. Still our regressor's RMS errors on such test data are no worse than the state of the art Artificial Neural Network regressors~\cite{lampeZiaukas2022lstmRoadFriction,hanyang2018annFriction}.
\end{abstract}

\section{INTRODUCTION}

Tracking the fluctuating friction between a road surface and passing vehicles is useful for two reasons. Firstly drivers can be alerted to potentially slippery road sections, and secondly winter road maintenance operations can be adapted to the actual 
surface conditions.

Narrow, curved and mountainous road sections present extra safety challenges, under wintry weather conditions.
A national road authority´s report~\cite{strandvik_haugvik_risikokurver_2018}
 documents a higher accident rate at such sections.
It shows that the risk of accidents increases, with any increase in the length and the curvature of the curve. Hence we need to track slippery conditions on curved roads. 

\subsection{The challenge of detecting slipperiness}
Modern vehicles come equipped with safety systems designed to reduce the risk of accidents occurring due to loss of traction. These systems include anti-lock braking system~(ABS), electronic stability control (ESC), and traction control system. The ESC feature has been mandatory in most cars and trucks, in many countries. 
The ESC system is activated when there is a danger of rollover or {\textit{dangerously excessive}} skidding. In the last two decades, ESC systems have led to a lowering of road fatalities. But these systems are {\textit{reactive,}} and trigger only after the vehicle has entered a dangerous situation, and the driver is about to lose control of the vehicle. Therefore  challenging winter conditions in difficult terrain demand solutions that are more predictive, and proactive.

\subsubsection{The vehicle state}
The ESC and ABS features use proprietary triggering rules, but basically use a common set of kinematic signals as inputs~\cite{rajamani2011vehicleDynamics}. These are vehicular signals that are affected by the state of the road surface, and we list them next.

\begin{figure}
{\begin{center}
\includegraphics[width=0.48\textwidth, keepaspectratio]{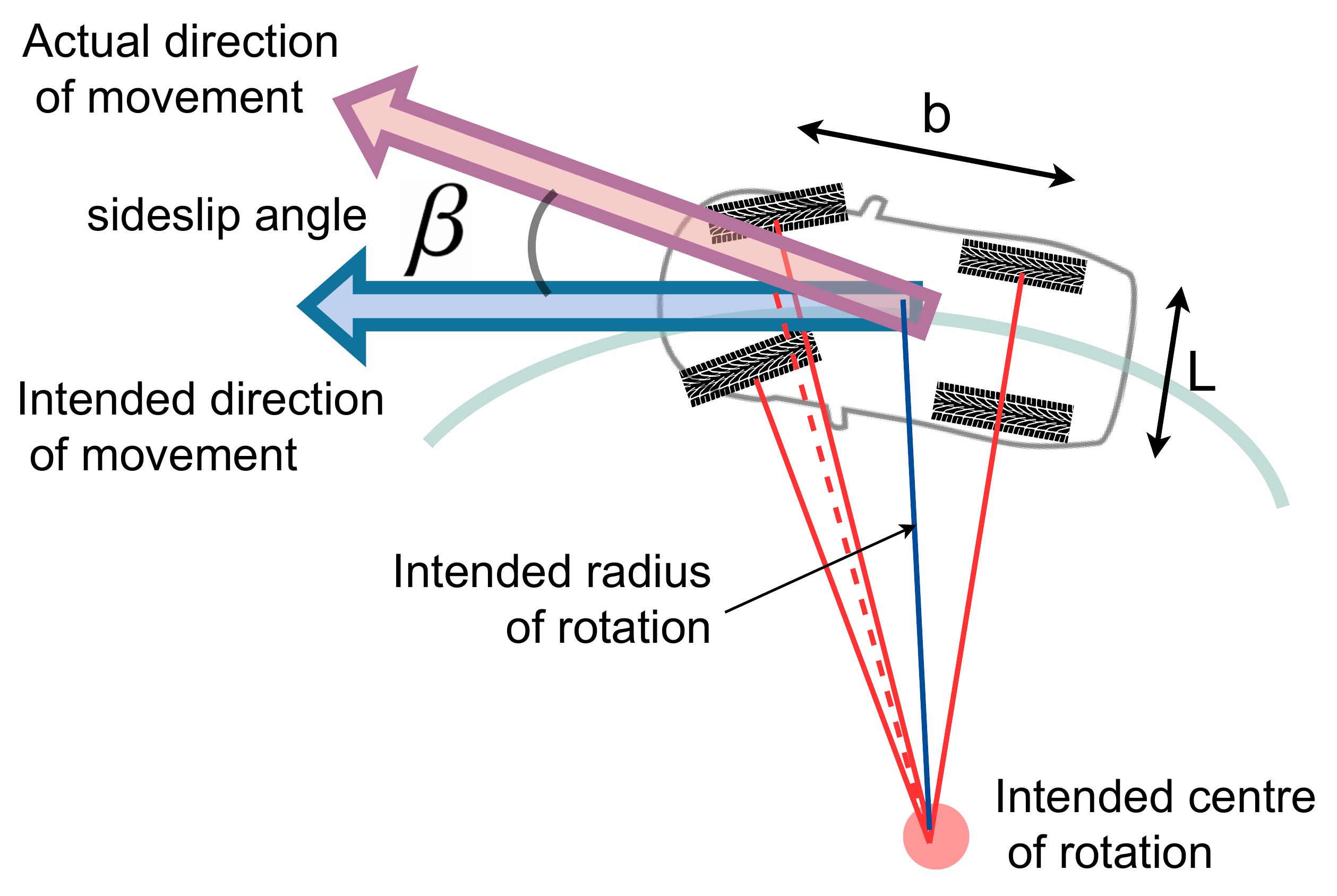}
\end{center}}
\caption{The sideslip, and the intended centre of rotation.}
\label{fig:sideslipIllustration}
\end{figure}

The slip of a wheel describes the movement of the tyre's contact patch relative to the ground. The {\textit{Longitudinal slip}} of a wheel is the difference between the speed of the wheel's axle, and the surface speed of the tyre due to wheel rotation alone. The {\textit{sideslip}} of a wheel is the angle between its intended direction of movement, and its actual direction of movement. Similarly the sideslip at any point of the chassis is the angle between the intended direction of movement at that point~(as determined by the intended longitudinal speed of the vehicle, and by the steering angle), and the actual direction of movement at that point in the chassis.
The intended yaw rate is the ratio of the longitudinal speed to the intended radius of rotation, which is shown in Figure~\ref{fig:sideslipIllustration}.
The {\textit{yaw rate}} excess is the difference between the intended yaw rate, and the actual yaw rate.

Thus the vehicle state is a real-time signal with fast changes, while the road state is a relatively slowly-varying parameter.

\subsubsection{The road state}
The {\textit{road state}} can be represented either by a discrete parameter to  classify the surface texture, or by the continuous parameter: the coefficient of friction offered to a  given tyre.
The longitudinal and lateral frictional forces have a nonlinear relationship with the longitudinal slip and sideslip respectively. The peaks of these 
relationship graphs determine the limits of adhesion, and these are directly proportional to the coefficient of friction. 
Hence it makes sense to track this coefficient. 

The coefficient of friction between the road surface and the tyre can vary with the vehicle type, the tyre type and also tyre wear.
This coefficient can also vary within a road section, depending on whether any water, snow, or black ice is spread uniformly over the road section or patchily. Nevertheless, at each designated road section,
and for each vehicle, we shall estimate the worst friction coefficient over the different patches within that section.

\subsubsection{Difficulty of estimating the coefficient of friction using existing sensors}
This friction parameter tracking task is hard because:
\begin{itemize}[align=parleft,left=0pt..1em]
\item 
as of now, there is no sensor that is cheap, 
and can directly measure the coefficient of friction,
\item
although the yaw rate excess can be easily derived from the IMU sensors, there is no cheap 
sensor available to {\textit{accurately measure}}  either the sideslip or the longitudinal slip. 
Indirect inference of these signals 
is difficult because:
\begin{itemize}
   \item
     when using an observer~(virtual sensor) that is based on the dynamics of vehicle motion~\cite{doumiati2012bookVehicleDynamicsEstimation}, the inference of kinematic signals requires knowledge of frictional forces on tyres. We have to either measure these forces or estimate them. Measuring them is neither easy nor cheap.  Estimating them requires reliable models of friction forces on tyres. These tyre-force models in turn require knowledge of the coefficient of friction.
    \item
    when using an observer that is based only on the kinematics of vehicle motion~\cite{farrellyWellstead1996lateralVelocityObserver,savaresi2017kinematicsBasedSideslip}, there is no need for models of frictional forces. This makes it easier, 
    but still the design of this observer has some limitations. Firstly current automotive grade IMU sensors  do not by themselves lead to estimates that maintain low error over time. And secondly these observers need tuning.
\end{itemize}
\end{itemize}
Therefore the problems of slip signals estimation and friction coefficient estimation have seen vigorous research efforts during the last two decades. Next we briefly review these developments.

\subsection{Previous works estimating the road state\label{section:previousWorksRoadState}}




Special sensors~\cite{tracSense,coventry2018virtualTyreForceSensors,brustad2020fieldStudyOfWinterRoadSensors,tongjiUniversity2022cameraAndVehicleDynamics,kth2021fusionOfHeterogeneousFrictionEstimates} have been developed for estimating road conditions, but these have been too expensive to be deployed in production vehicles today. Moreover, in the case of camera based sensors, we get only a limited
 accuracy~\cite{casselgren2021floatingCarData} in the face of variable lighting conditions.


\subsubsection{Dynamical models for indirectly inferring the road state}

Braking manoeuvres can be used to estimate some points on the graph of the friction force versus longitudinal slip characteristic~\cite{gustafsson1997slipBasedTireRoadFrictionEstimation}. But this requires estimating the braking forces accurately and the estimators/observers to be tuned to the vehicle parameters.

Parametric models for tyre forces can be incorporated into observers for the sideslip and other vehicle states, and ultimately the friction state of the road as well~\cite{antonioLoria2018switchedObseerverForFriction,baffetCharara2009estimationOfWheelSideslip,doumiati2012bookVehicleDynamicsEstimation,gripFoss2008automaticaNonlinearSideslipEstimation}.

A commercial service that estimates the friction coefficient experienced by a vehicle, is offered by 
the company NIRA dynamics~AB~\cite{niraDynamics,casselgren2021floatingCarData}.
This service is 
 based on gathering data from existing vehicular sensors, 
 and communicating them to a cloud server. Their algorithm has two components: an onboard part and  a cloud server part. The onboard part processes  vehicular data gathered over an Onboard diagnostics~(OBD) interface, from  those time durations when the vehicle is accelerating or braking. 

\subsection{Previous works estimating the vehicle state}
Special sensors for the use of GPS/GNSS modules have been proposed~\cite{berntorp2016wheelSlipEstimationInertialGPSwheelSpeed,mazzili2021smartTyreTechnology,
ucl2018frictionalForceAndWheelSlip,viehweger2021vehicleStateAndTyreForceEstimationDemonstrationsAndGuidelines,changchunWaterloo2021smartTyre} for the estimation of the vehicle state. But
 these come with costs and limitations.

The vehicle sideslip is the most challenging of our vehicle states to estimate. This estimation in turn requires the estimation of the longitudinal and lateral velocity signals.
For this some researchers~\cite{doumiati2012bookVehicleDynamicsEstimation,haudum2018vehicleSideslipAndBank,borelli2019adaptiveSideslipEstimation} have constructed observers, based on
dynamical models of vehicle motion,
requiring models of tyre-road friction forces.

More promising for us are observers constructed based on kinematics models of vehicle motion~\cite{marco2020imuPlusLevelSensorsVehicleMotion,savaresi2017kinematicsBasedSideslip}. The observer of~\cite{savaresi2017kinematicsBasedSideslip}  only requires standard sensors that are available in production vehicles.

\subsection{Problem formulation and our contributions}
The problem is to design a system to collect location-specific, vehicle state data from several vehicles, and to calculate in nearly real-time, estimates of the slipperiness of chosen road sections. The constraints on the system are the following: (i)~vehicle states should be estimated onboard the vehicles, in nearly real-time,   and using only sensors that are available in vehicles in production today, (ii)~the part of the system that is onboard vehicles should be of a purely monitoring nature, and shall not interfere with the drivers'  inputs - there can be no demand for braking manoeuvres, and (iii)~communications between vehicles and the infrastructure should happen over an existing V2X wireless protocol such as WiFi~(ITS-G5, or DSRC), or cellular communications~(4G/5G).

\subsubsection{Our approach}
At curved road sections, the sideslip and yaw rate excess signals seem to be the key kinematic signals affected by friction. If we fix the speed of the vehicle, then the amplitudes of these signals vary inversely with the friction coefficient.
 Hence our premise is that:  {\textit{even some  
 reasonably accurate estimates of these signals
  shall be strongly correlated with the friction coefficient.}}

Suppose that we were to run an observer~(virtual sensor)
for these signals, only over the short duration that  it takes to traverse a given curved section.
Then we can reduce inaccuracies due to gyroscope drift in the IMUs. And this could help
to build a reasonably accurate estimate of the sideslip signal at curved sections.



\subsubsection{Our contributions\label{section:ourApproachAndContributions}}
We make two contributions.

The first contribution consists of 
    our system architecture for collaborative friction estimation, to deliver estimates in nearly real-time. Its merits are that: we shall only use existing vehicular sensors, and our system
        generates only a small amount of data exchange over the V2X communication links. In specific, when a  vehicle passes through a designated road section, the vehicle and the infrastructure exchange at most twenty standard-sized WiFi data packets.

The second contribution is the set of
    our findings through simulations.
 We find that our system  estimates the friction coefficient with a low magnitude of worst case error, which is less than 10\%, and an even lower average value for the magnitude of the error.
Another finding shows that our system
    can be implemented and maintained with low complexity and effort.
  In specific, there is only a marginal drop in accuracy if we deploy  at all curved road sections, a common global regressor for the friction coefficient. This choice increases the size of the training set. A local regressor can only be trained on ground truth  friction data collected at that specific road section. But the global regressor can be trained  the combined  ground truth data from all the road sections.


%

\begin{figure*}
\begin{center}
\includegraphics[width=0.99\textwidth]{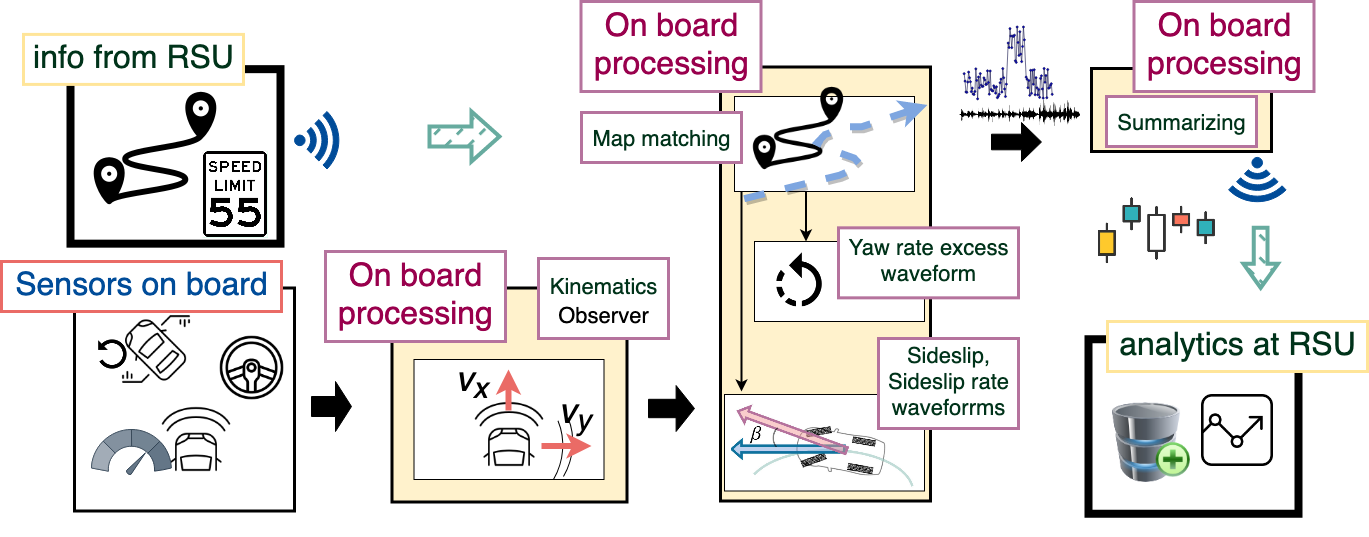}
\end{center}
    \caption{The data flows in our system architecture}
  \label{fig:theCriSpStack}
\end{figure*}

\section{Our system architecture}
The architectural form and the algorithmic structure of our system follow from its constraints and  function.

At the design stage, we must list a  set of curved road sections, where we shall deploy our system. Straight road sections are excluded.  This excludes the practically important case of  straight bridges that can get icy under wintry weather.

\subsection{Inputs and outputs}
Our system is a distributed one. It has processing and communicating modules onboard vehicles that are enrolled in it. And it has  processing and communicating modules at the infrastructure end. We shall assume that the infrastructure nodes take the form of {\textit{Road side units}}~(RSUs) that are installed at the side of designated road sections.
But our system can function essentially the same way if the infrastructure nodes are all merged into a single cloud server, that is accessible over a cellular communication network~(see Section~\ref{section:wirelessStandardDesignChoice}). Even when the cellular option is what is deployed, the following discussion can proceed as if there is a virtual RSU at each designated road section, delivering the same outputs that a physically installed  RSU would.

When our system is configured and  deployed, the real-time inputs to our system are the vehicular sensor measurements. The real-time outputs are {\textit{intervals}} of friction coefficients computed at the RSUs. Each designated road section has a dedicated RSU that tracks a interval of friction coefficients at that section.
\subsubsection{interval of friction coefficient values}\label{section:intervalsOfFrictionConditions}
Consider some fixed surface condition, at a given road section.  
Then the friction coefficient between the surface and any tyre depends on factors such as: the vehicle type, the tyre type, tyre pressure, and the level of  tyre wear. The last of these factors is much more significant than the others.
 There is no significant difference in  friction coefficients experienced by new issues, from the majority of different tyre types that are used in a given climatic season\footnote{In some parts of the world, a new issue of a  typical tyre used in winter can offer a noticeably higher friction coefficient that a new issue of a typical tyre used in summer.}.
 A badly worn-out tyre experiences a friction coefficient that is about eighty percent of that experienced by a newly bought issue of the same tyre type~\cite{wright2019tyreAgeAndWear}. Similarly, a drop in pressure from~2.5~bar to 2~bar causes the friction coefficient to increase by about eight percent~\cite{braghin2006tyreWearModel}.

 Therefore we can parametrically describe the interval of typical friction values as the interval:
 \begin{gather*}
     \left[  \,
        0.8 \times \mu_{\text{new}}\left( t \right) , 
        \, \mu_{\text{new}}\left( t \right)
        \,
     \right] ,
 \end{gather*}
where~$ \mu_{\text{new}}\left( t \right)  $ is the friction coefficient experienced at time~$ t , $ at the given road section,
by a new issue of a typical tyre type, and at a pressure that is slightly lower than that recommended by the tyre manufacturer. This way, we have a reasonable model to capture the variations of friction coefficient, due to variations in both wear level and tyre pressure.

We shall now give an overview of the data transformations and flows in our system, which are  illustrated in Figure~\ref{fig:theCriSpStack}.

\begin{figure}[!htb]
\begin{center}
  \includegraphics[trim={6cm 4.5cm 0cm 0cm},clip, width=0.23\textwidth]{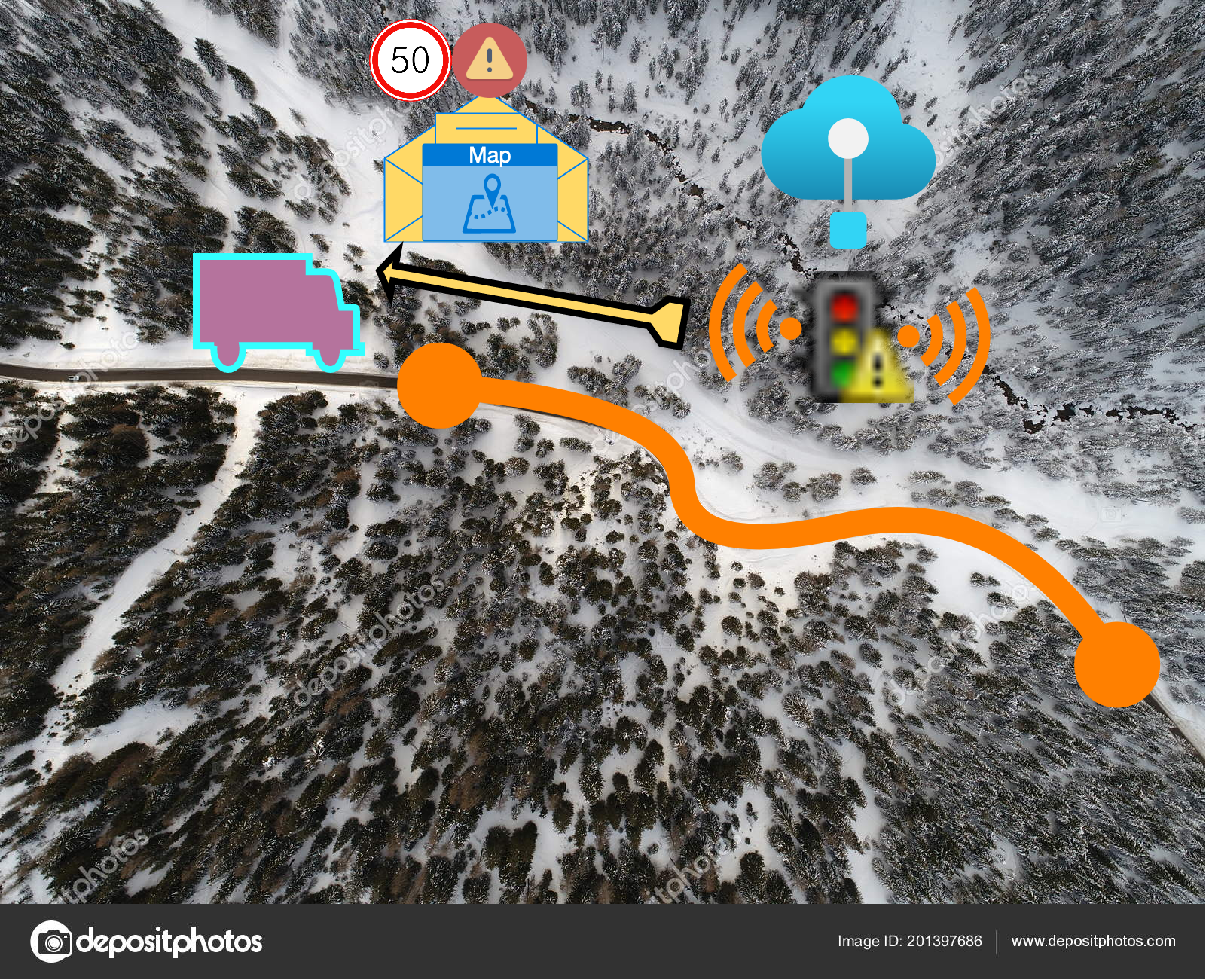}
  \includegraphics[trim={6cm 4.5cm 0cm 0cm},clip, width=0.23\textwidth]{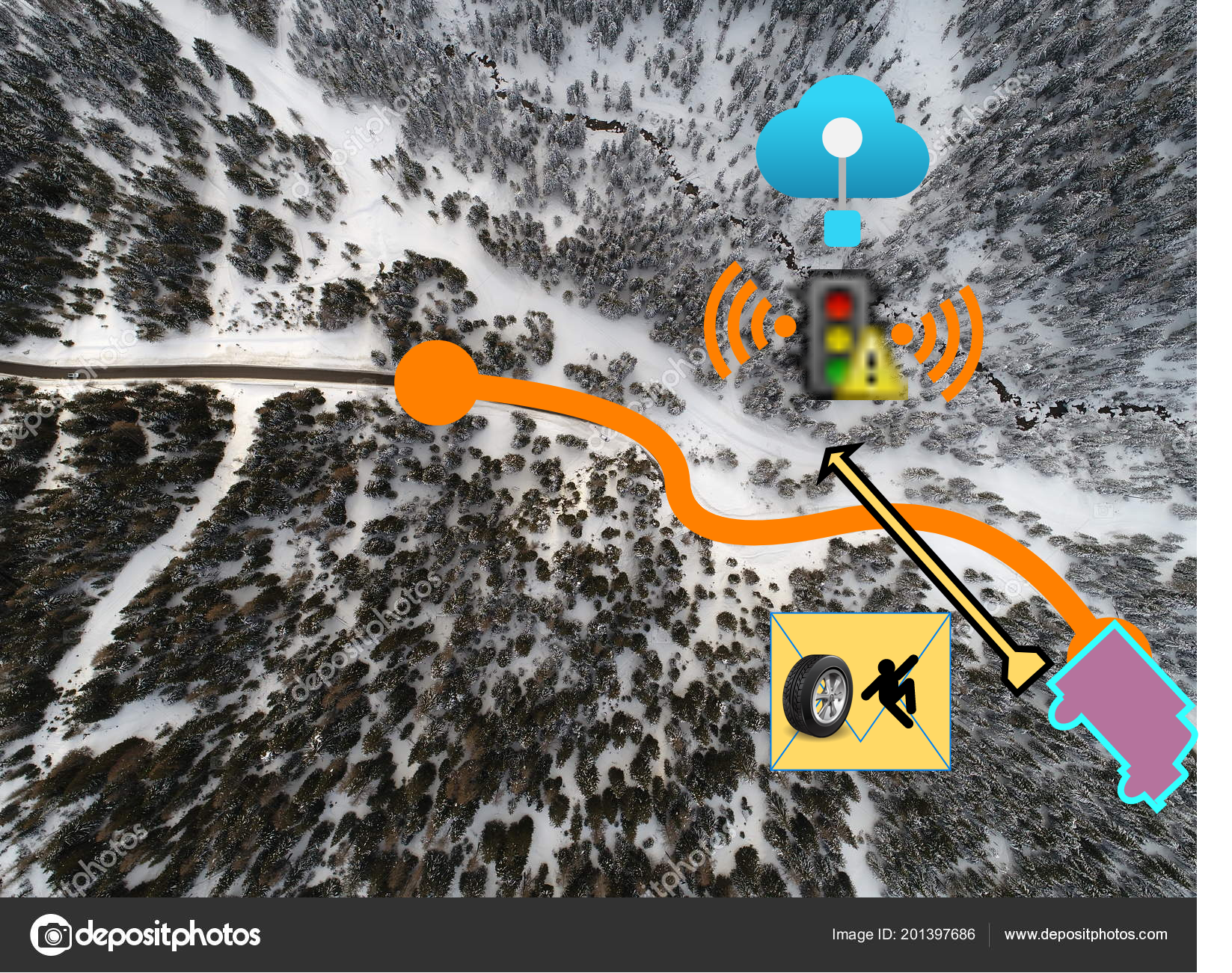}
\end{center}
    \caption{The vehicle receives a message from the RSU upon entry to the road section, and broadcasts a message to the RSU when leving the section.}
  \label{fig:vehicleRSUMessaging}
\end{figure}

Our system functions as follows. A new arrival of a vehicle at a designated road section is detected either by a GPS-linked module on the vehicle itself, or by a road surface sensor such as a magnetic loop. This detection triggers a message from the RSU to the vehicle~(see Section~\ref{section:messageFromRSUtoVehicle}). This message in turn triggers and seeds some processing onboard the vehicle. The end result of that processing is a message back to the RSU, with a {\textit{kinematic summary,}} which is a highly compressed digest of the vehicle state signals~(see Figure \ref{fig:vehicleRSUMessaging} and Section~\ref{section:messageFromVehicletoRSU}). The RSU then uses the received kinematic summary as either the whole or part\footnote{In those design variations where the kinematic summary is not the whole of the feature vector, the rest of it is a summary description of the road geometry.} of a feature vector for a regressor, that shall estimate the friction coefficient experienced
by the individual vehicle. 
Finally, the RSU applies its rule for combining
 friction coefficient estimates for recently passed vehicles, to estimate the road section's present interval of friction coefficients.

\subsubsection{Message from RSU to arriving vehicles\label{section:messageFromRSUtoVehicle}}
  The content of this message is independent of the type of the vehicle receiving the message.  The message comprises:
\begin{description}
%
    \item[$\blacklozenge$ ] {\textbf{speed advisory:}} maximum allowed or recommended speed, which depends on the RSU's estimate of the 
    road conditions.
    \item[{$\blacklozenge$} ]{ {\textbf{map:}} of the 3D geometry of the road segment. 
    The lanes are to be specified as geometric curves in space~(such as a patchwork of splines), and possibly augmented with GPS coordinates of some points on the road. The length of the road segment should of course be included.}
\end{description}
        
\subsubsection{Message from leaving vehicles to RSU \label{section:messageFromVehicletoRSU}}
 This message is also independent of the type of vehicle sending it. It depends solely on the kinematic signals of the vehicle at this road section. It conveys:
\begin{description}
  \item[{$\blacklozenge$}]{{\textbf{speed taken:}} The profile of speeds taken at this road section. In specific, it shall provide a list of quantiles of the speed waveform over the section. The waveform we shall use shall actually be a filtered estimated derived from the four wheel encoders, as done in Selmanaj~et~al.~\cite{savaresi2017kinematicsBasedSideslip}.}
    \item[{$\blacklozenge$}]{{\textbf{vehicle state:}} some chosen set of quantities that represent the 
    skidding that was experienced. What exactly should be computed and packaged as the Kinematic summary? This is the main design question concerning the content of this message.} 
\end{description}

Section~\ref{section:onBoardProcessing} lists the processing stages that take place onboard the vehicle
between the reception and sending of the above two
messages. Next we see the main algorithmic pieces in our system.



\subsection{Overall algorithmic structure}
Those design choices that are locked-in have been guided by the following three considerations: (i)~we shall minimize the amount of data exchange between vehicles and RSUs, (ii)~we shall make use of models of vehicle motion, only for estimating the vehicle state signals, for which the models can provide reasonable accuracy,  and (ii)~we shall make use of data-driven machine learning algorithms only for learning the road state, for which vehicle motion models cannot give accurate estimates given our constraints. The last consideration arises because: (a)~we are constrained to use only data from: IMUs, steering angle sensors, and wheel speed encoders, and (b) it is very hard to acquire accurate models for tyre-road forces, or to accurately estimate these forces. Because of the last constraint, we abandon dynamics models of vehicle motion, and instead use kinematics models to estimate the vehicle state.

Our system architecture splits the overall task of estimating friction coefficient intervals, into three sub-problems, listed below.
\paragraph{Sub-problem~1: Model-driven vehicle state estimation}
{The kinematic summaries shall be primarily based on the yaw rate excess  and the sideslip angle signals. The yaw rate excess is easy to estimate, but the sideslip  is not. When a vehicle takes a curved trajectory, then the sideslip signal is observable from the IMU signals.
We shall use the observer of  Selmanaj~et~al.~\cite{savaresi2017kinematicsBasedSideslip}, that is based on the kinematics of vehicle motion. It requires only data from sensors that we have chosen, and has an acceptable accuracy which has been validated in actual road tests.
}

\paragraph{Sub-problem~2: Data-driven Regression of friction coefficients}
{We shall infer the coefficient of friction experienced by an individual vehicle, using a machine learning regressor. The exact contents 
of the feature vector is not fixed apriori. We shall explore using either just  the kinematic summary as the feature vector, or this summary augmented by information about road geometry. Training the regressor
requires ground truth data from all designated road sections. In specific, we require the road authority to make several hundreds of runs of vehicles, under different friction conditions,  where the friction coefficients in operation are accurately measured with specialized equipment~(see for example~\cite{viaFriction}).
}

\paragraph{Sub-problem~3: Collaborative estimate at the RSU}{The collaborative estimate should be a function of the received sequence of individual estimates, of friction coefficients experienced by vehicles at that section.
We choose one of the simplest possible rules for estimating the current interval of friction coefficients.
The RSU at a given road section shall simply compute the median of the recent sequence of individual friction coefficient values.
For this we choose a size for the time window, an hour for example.
Then the resulting average number of vehicle passages per hour, and the accuracy of the interval estimate depends on the usual intensity of traffic at that road section\footnote{during periods of light traffic, for example in the middle of the night, the accuracy of the RSU's interval estimate may suffer, because of lesser traffic intensity.}.
We make three assumptions in this context: (i)~the friction coefficient interval may only change very little over short durations such as an hour, and (ii)~on average, at least a few tens of vehicles pass through any designated road section, (iii)~the variation in tyre wear and tyre pressure among the vehicles is such that, their friction coefficients at any road section have a probability distribution 
that is symmetric about the midpoint of the interval. The last assumption implies that the RSU's interval estimate is completely fixed, once the median estimate is fixed.
}

\subsection{Processing stages onboard vehicles\label{section:onBoardProcessing}}
The main tasks of onboard processing are to: (i)~isolate and extract sensor measurement waveform segments, from only the time duration when the vehicle is located within the spatial bounds of the designated road section, (ii)~estimate from these waveform segments, the vehicle state signals, and (iii)~compress the vehicle state signal waveforms into the kinematic summary. These tasks are achieved by the following processing stages.
\paragraph{Stage 1: Trajectory estimation and Map matching}

In GPS-enabled vehicles, this stage reduces to map matching alone, because the GPS module readily gives a good enough estimate of the trajectory. 
We simply compare the GPS coordinate of the vehicle's trajectory with the start and end points of the road section, as given on the map sent by the RSU.
Once the end point has been crossed, the subsequent stages 
can be triggered.

On the other hand, in vehicles without GPS, or in locations where the GPS signal reception is poor, this stage is more involved.
This stage has a component that runs continuously all the time, 
estimating the spatial trajectory, by IMU-supported dead reckoning. Essentially, we shall estimate the spatial coordinates of the trajectory, by integrating the
kinematic observer of Selmanaj~et~al.~\cite{savaresi2017kinematicsBasedSideslip} for the longitudinal, and lateral velocities, with an easy augmentation of this observer with an element for the vertical velocity.

With this trajectory estimate, we shall do map matching
by identifying special way points on the map. For example,
we can identify points of maximum curvature, the signed value of this curvature, and the distances of these points from each other and the start and end of the road segment. Then the map matching problem reduces to matching  waypoints on the map~\cite{icmlt2022roadQualityMeasurementIMU}.

\paragraph{Stage 2: Estimating the vehicle state}
The sideslip and the yaw rate excess are
the vehicle signals that are most affected by friction at curved sections\footnote{While braking, the longitudinal slip and longitudinal acceleration are also affected.}. Indirectly, as a consequence of the driver's response to the curve, the longitudinal speed is also affected by the friction. Another relevant signal is the time derivative of sideslip, which is a byproduct of the observer used for estimating the sideslip.
The above mentioned four kinematic signals shall comprise our candidate vector for the vehicle state.

The longitudinal velocity can be estimated from the wheel speed sensors. 
The yaw rate excess can be estimated using this  velocity, the steering angle, and the IMU output. The sideslip requires both longitudinal and lateral velocities. The details are given in 
Section~\ref{section:observer}

\paragraph{Stage 3: Summarizing the vehicle state signals}
Our aim here is to compress the waveforms at the output of the previous stage. These discrete-time waveforms have a sample rate of 100~Hertz.
We shall strip away the timestamps, and form a bag of sample values from each of the four scalar waveforms. We shall then represent each waveform by a statistical summary of its bag.
A minor design choice concerns what the statistical summary should be. We shall explore in our simulation studies the possible formats: (i)~mean and standard deviation, (ii)~quantiles, and (iii)~skewness and kurtosis.

\subsection{Regressor: local or global?\label{section:designChoices}%
}
The main design choice we face arises in Sub-problem~2.  This choice determines the overall formulation  of the learning problem.

Even between locations that are just a few hundred metres away from each other, surface conditions can be different, because sunlight, road wetness, temperature etc. can be different.

Should  each individual road section have a {\textit{local}} friction coefficient regressor, 
  that is trained on ground truth data from only that road section, 
  or can we use a {\textit{global}} regressor that covers all the possible road sections?

Under the first alternative, the road operator has to train  a collection of regressors, each specific to a corresponding individual road section. 
For every regressor, the feature vectors only carry a summary of the kinematic signals, and contains no explicit information about the particular road section. Therefore, the training data for a regressor has to come from vehicular runs on the specific road section for which the regressor is being trained. This specificity reduces the size of the training set available for each regressor. 

Under the other alternative, the road operator has to train a single  regressor that can predict the road state at all road sections. 
But then the 
feature vector has to include entries that capture the geometry of the road. Under this choice, the regressor can be trained on one large corpus of data from vehicular runs on diverse, representative road sections. The potential advantages are that: (a)~we may get a higher accuracy of the regressor because of the enlarged size of the combined training corpus,
and (b)~the road operator can deploy an identical copy of the trained regressor at each road section.


\section{The observer for the vehicle state\label{section:observer}}
\subsection{Observer for longitudinal and lateral speeds}
We 
use the kinematics-based vehicle state observer 
of Selmanaj et al.~\cite{savaresi2017kinematicsBasedSideslip}. 
The observer takes as its inputs measurements from the IMU, 
and wheel speed encoders. The observer produces an estimate of the longitudinal and lateral vehicle speeds. Below is the evolution equation, where the orange coloured terms highlight the synthetic changes made by Selmanaj et al.~\cite{savaresi2017kinematicsBasedSideslip} to the original Kinematic laws of motion:
\begin{gather*}\label{eq:kinematicObserver}
\begin{split} 
\begin{bmatrix}
    \dot{\hat{V}}_x(t) \\ \dot{\hat{V}}_y(t) 
\end{bmatrix} = &
\begin{bmatrix}
    {\textcolor{orange}{-\alpha_0-\alpha_1\left\lvert\omega_z(t) \right\rvert }} & \omega_z(t)  \\
    - \left( {\textcolor{orange}{\alpha_2}} + 1 \right)\omega_z(t) & \textcolor{orange}{-F(t)}
\end{bmatrix}
\begin{bmatrix}
    \hat{V}_x(t) \\ \hat{V}_y(t) 
\end{bmatrix}\\& +
\begin{bmatrix}
    1 & 0 \\ 0 & 1
\end{bmatrix}
\begin{bmatrix}
    A_x(t) \\ A_y(t) 
\end{bmatrix}
+
\begin{bmatrix}
    {\textcolor{orange}{\alpha_0+\alpha_1 \left\lvert \omega_z(t) \right\rvert }}
    \\ 
    {\textcolor{orange}{\alpha_2 \omega_z(t) }}
\end{bmatrix}
V_x^{\text{encoder}} .
\end{split}
\end{gather*}
Here~$\omega_z$ denotes the yaw rate in the road frame. Similarly, $A_x$ and $A_y$
 denote longitudinal and lateral acceleration 
 in the road frame. These are obtained from the IMU 
 measurements, after correction for the chassis roll angle and the banking angle of the road. We assume that the banking angle profile has been supplied by the RSU. We estimate the chassis roll angle by applying a complementary filter on the IMU.
 
A longitudinal speed estimation block~\cite{savaresi2017kinematicsBasedSideslip} estimates the vehicle speed $V_x^{\text{encoder}}$ using the wheel speed encoder measurements, the steering angle, and yaw rate measurements. 
The coloured parts of the observer's RHS are correction terms that stabilize it during straight driving, when they push the lateral speed towards zero. 

Our implementation has  no correction for offsets
between the location of the IMU and the vehicle's centre of gravity.
This is because our simulated IMU is placed right at the centre of gravity of the vehicle. 
The observer constants were tuned to the simulated vehicle, by minimizing the root mean square~(RMS) error of the sideslip angle estimate. See Table~\ref{table:tunedObserverConstants} in the appendix for a list  of the observer parameters used. 

\subsubsection{Smoothing using a two-way filter}
Onboard the vehicles, the sensor signal record over the road section has to be isolated first before processing (see Section~\ref{section:onBoardProcessing}). 
Hence the signal processing operations performed on the sensor signals do not have to be causal. In specific, we can perform smoothed estimates instead of filtered estimates. This means that our low pass filters, high pass filters and complementary filters can be run twice: once in the forward direction of time, and again in the backward direction of time. This two-way filtering eliminates the phase lag in estimates.

We apply this two-way filtering to remove noise without adding phase lags, in estimating the following signals:
$A_x , A_y$ (longitudinal and lateral acceleration),
$ \delta $ (steering angle),
$ V_x^{\text{encoder}} $    (the longitudinal speed estimate from the longitudinal estimation block),
 sideslip rate estimate~\eqref{eqn:sideslipRate}.
and the estimated longitudinal velocity and lateral velocity of the observer.

A two-way low-pass filter with 100Hz sampling frequency, 10Hz cutoff frequency and a filter order of 30 was used in all of these estimations. 
For chassis roll angle estimation, the low-pass filter is used in the complementary filter. 

\subsection{Estimating yaw rate excess}
The yaw rate excess is the difference between the  measured yaw rate of the vehicle, and the intended yaw rate.
The intended yaw rate is calculated from the steering angle and the longitudinal speed, based on the geometry shown in Figure~\ref{fig:sideslipIllustration} as:
\begin{gather*}
   \frac{2V_x^{\text{encoder}} \ \tan(\delta)}{2b+L \ \tan(\delta)} ,
\end{gather*}
where 
$\delta$ is the wheel steering angle~(the average of steering angles of the front wheels), $L$  the wheel track, and $b$  the wheel base. 

\subsection{Estimating sideslip and its rate}
\begin{figure}
{\includegraphics[width=0.40\textwidth, keepaspectratio]{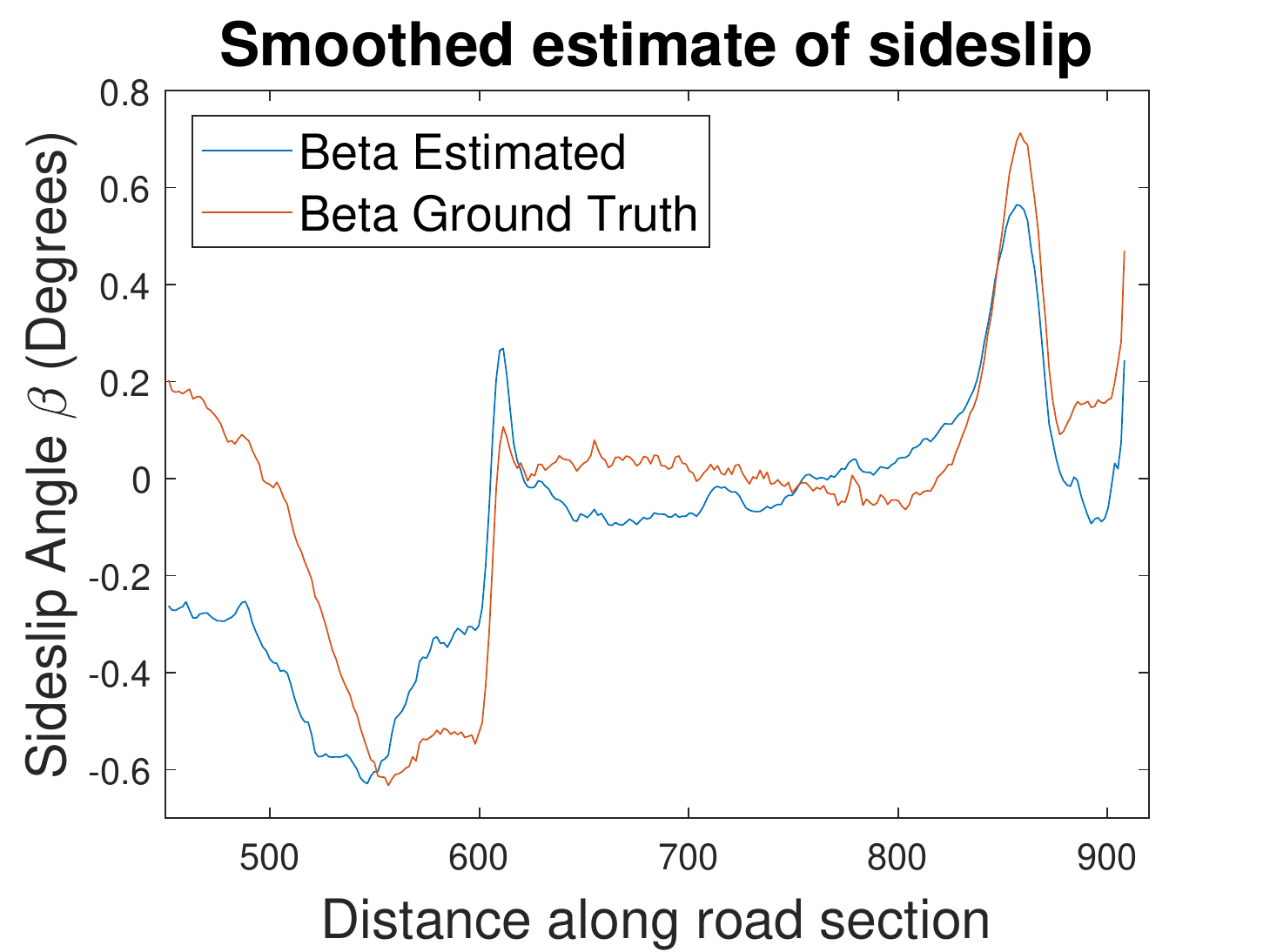}}
\caption{Two-way filter smoothed estimate of sideslip signal, from a simulated run on the S-turn scenario with uniform road friction coefficient of 0.2 and 72kmph vehicle speed.}
\label{fig:sideslipEstimate}
\end{figure}

The sideslip estimate is:
\begin{align}
    {\widehat{\beta}(t)}
    & = \arctan{
       \left( \frac{{\widehat{V}}_y(t)}{{\widehat{V}}_x(t)} \right) 
       } .
\end{align}
Figure \ref{fig:sideslipEstimate} shows the sideslip waveform of a simulated vehicle. The performance is 
acceptable,
as the peaks are picked up pretty well. An estimate of the rate of sideslip angle 
is calculated with:
\begin{align}
\label{eqn:sideslipRate}
   {\widehat{\dot{\beta}}(t)} &  = 
   {\frac{ A_y(t) - \omega_z(t) V_x^{\text{encoder}}(t)}
   {V_x^{\text{encoder}}(t)}
   } .
\end{align}

\section{Generating training data from simulations}\label{section:simulationDescription}
At the outset, four vehicle simulators were considered: Carla, CarMaker, Dyna4\cite{noauthor_dyna4_nodate} and the Vehicle Dynamics tool box available in Simulink. Carla appeared to have a primary focus on autonomous driving with related sensors. Simulink's Vehicle Dynamics toolbox was found to lack certain assistive driving systems such as ESP and ABS. CarMaker was too expensive, 
 for our budget. 
 We selected Dyna4, as it provided vehicle dynamics models and the vehicle sensors that we required at an affordable price. 

\subsection{Driving scenarios from a road database}
The Norwegian nasjonal vegbank \cite{noauthor_nasjonal_nodate} is a public road network database, and Vegkart \cite{noauthor_vegkart_nodate} is a companion tool for visualising and accessing the data.
We used the report~\cite{strandvik_haugvik_risikokurver_2018} on accident risks at curved road sections, 
as a starting point for scenario selection. The report states that as the curve length increases from 50m to 200m, and the radius of curvature reduces from 500m to 50m, there is a significant increase in the number of accidents. Vegkart was used to select three curved road sections that fitted these criteria. 
The rated speed of the curved scenarios were calculated from Equation \ref{eq:bankAngle} in the appendix. See section \ref{section:roadQueries} in the appendix for links to the selected road section queries in Vegkart. The scenarios are:
\begin{itemize}[align=parleft,left=0pt..1em]
    \item \textbf{Long turn}: contains a 600m long curved section, in which the turning portions have a radius of curvature between 100-150m. The rated speed is 80kmph. 
    \item \textbf{S-turn}: contains two curves where each curve is  close to 300m long, and the turning portions have a radius of curvature between 55-187m. The rated speed 
    is 70kmph.
    \item \textbf{Sharp turn}: contains a 100m long sharp curve, where the radius of curvature is between 20-43m. The rated speed is 20kmph. 
\end{itemize}

\subsubsection{Scenario reconstruction from database}\label{section:scenarioReconstruction}
The reconstruction of the road scenarios in the Dyna4 simulation environment took several steps. First the data from the selected road sections were exported. The exported data includes the three spatial coordinates of the road lane's centre, the radius of curvature, and the speed limit for each road segment. The road bank angle was not available in NVDB. To remedy this, 
we calculated the bank angle 
using Equation~\eqref{eq:bankAngle} in the appendix,
which follows the road authority's design guidelines~\cite{ziani_overhoyde_nodate}. The bank angles were capped to 8~degrees, as that is the maximum 
value allowed in the design guidelines.


\subsection{Simulation configuration and runs}
Recall that the overall architecture of our system is shown in  Figure~\ref{fig:theCriSpStack}. In our simulations, 
two parts of the architecture are not considered.

Firstly we abstract and idealize the V2X communications.
Secondly we assume that the map~matching step is executed
with negligible error, because good validated algorithms exist~\cite{icmlt2022roadQualityMeasurementIMU}.

A pre-configured Volvo SUV XC90 2015 vehicle model from Dyna4 is driven through each scenario. The vehicle tyre model is the TMEasy 5 tyre model which allows for precise contact forces and simulation of inflation pressure. Noise is introduced to the data similar to that seen when driving on good asphalt. A vehicle driver model configured to be representative of skilled drivers is tasked with following the road lane. 

Each scenario was simulated following a test matrix consisting of three parameters: base value of friction coefficient, vehicle speed, and tyre wear. The base friction coefficient ranged from 0.2 to 0.7, with increments of~0.05. For each value of the base friction coefficient,  1000 vehicle runs were simulated. The vehicle speed was sampled from a triangular distribution centred on the rated speed, 
with lower and higher limits that are 20\% away from the rated speed. As mentioned in Section~\ref{section:intervalsOfFrictionConditions}, we model the effect of tyre wear and tyre pressure on the coefficient of friction by multiplying the base friction coefficient with a number uniformly sampled between 0.8 and 1. The driving scenarios were simulated with constant road surface friction throughout the section. In addition, the S-turn scenario was simulated one more time, but with lower friction in one turn and a higher friction in the other turn. 


Kinematic summaries were collected from a total of 44,000 simulated vehicle runs on four different scenarios. Each kinematic summary consist of the following signals: The wheel steering angle, the vehicle sideslip angle estimate, the sideslip angle derivative estimate, the yaw excess and the vehicle speed. Each signal is statistically summarised with the mean, standard deviation, the 20th, 40th, 60th and 80th quantiles, min, median, max, skewness and kurtosis.
Every kinematic summary was labelled with the lowest experienced friction coefficient as the ground truth. Each simulated scenario had a dataset of 11.000 vehicle runs which was shuffled and split into a training set (90\%) and a test set (10\%). In addition, a combined dataset was created from all of the scenario training and test data splits. The combined training dataset was shuffled after merging to ensure an even distribution of scenario data when performing cross-validation.

Next we describe the results of using the simulation setup to evaluate our design variations.
\section{Training the ML estimators}
Our solutions to sub-problems~2 and~3 depend on how accurately we are learning  to compute the friction
coefficient. 
The main design choices in our solutions are about the content of the feature vectors. 
Here we describe our investigative questions, and the answers to these from our simulation studies.


\subsection{Machine learning algorithm selection}
A set of ML algorithms were selected which consisted of Support Vector Machine~(SVM), Multi-layer Perceptron~(MLP), Random Forests, Gradient Boosting and XGBoost (Extreme Gradient Boosting) with Random Forests as the underlying tree method. All algorithms except for MLP  and  SVM were ensemble methods. The latter are preferred, as they generally perform better, at the expense of additional computation. In our case, the trained regressor shall be deployed on dedicated hardware, therefore the extra computational burden is acceptable, in exchange for better accuracy.

The ML framework SKLearn \cite{scikit-learn} and the gradient boosting library XGBoost~\cite{Chen:2016:XST:2939672.2939785} were used to perform the ML study. Plotly \cite{plotly} was used to generate the graphical plots. 
Each algorithm was configured to use the same values of hyper parameters where possible, namely: 100~estimators and a max tree depth of~6. Eight-fold cross validation was performed on the combined training dataset. Mean squared error was used as the loss function, and the following metric used for  error comparison between algorithms:
\begin{equation}\label{eq:MAPE}
{\text{MAPE}} = \frac{1}{n}\sum^n_{i=1}\left|\frac{y_i-\hat{y_i}}{y_i}\right|
\end{equation}
where $n$ is the number of samples, and the real numbers $y_i, \hat{y_i}$  are the true  
and 
estimated values respectively.

Figure~\ref{fig:modelCVPerformance} in the appendix shows a comparison of the algorithms' performances. The best performing algorithm was XGBoost (median MAPE score of 1.68\%), and the worst performing was MLP (median MAPE score of 2.84\%).
Therefore XGBoost was selected to be the 
algorithm of the road state regressor.

\subsection{Sequential feature selection}
The dataset contains 82 input features in all. 
To reduce the complexity of the resulting regressor, a sequential feature selection (SFS) algorithm was used. SFS starts with an empty set of input features, and greedily adds features that give the best score on the validation data. 
SFS stops when the desired number of features is reached. The SFS algorithm was run eight times with different numbers of desired features. The algorithm was used on two variations of the dataset. One where road geometry features were included, and one without. Results showed that the same kinematic features were selected for both datasets up until feature number 10. After that the selected features for each dataset diverged, and road features started to get selected - see Figure~\ref{fig:selectedFeatures} in the appendix. 

The performance of the regressor using the  growing sets of selected features was tested through 8-fold cross validation
 - see Figure~\ref{fig:nmbFeaturePerformance}. The smallest feature set of 3~features have a median MAPE close to 5\% while the largest set of 30~features is close to 2.5\%. After the selection of about 15~features, there are diminishing returns for both dataset variations. The set of 15 selected features was used for the remainder of the ML study as it provides good performance, the MAPE being close to 2.5\%, while also keeping the complexity of the regressor low. 

\subsection{Variations in statistical summaries}
\begin{figure*}
{\includegraphics[width=0.90\textwidth, keepaspectratio]{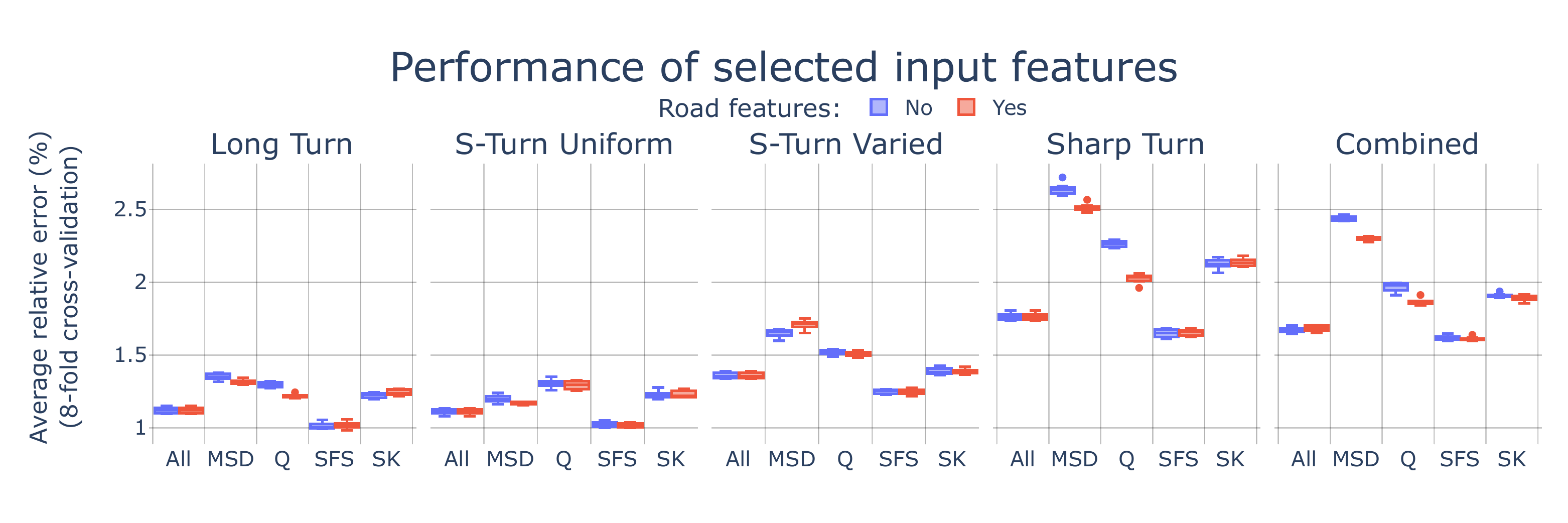}}
\caption{Performance of the XGBoost regressor, for variations of  
statistical summary formats in input feature vectors. The variations from left to right are: All features, mean and standard deviation, quantiles, SFS selected features, and skewness and kurtosis.}
\label{fig:inputFeatureComparison}
\end{figure*}

In addition to the features selected through SFS, a set of other input feature variations were examined: (i)~mean and standard deviation, (ii)~quantiles, and (iii)~skewness and kurtosis. For baseline comparison a feature set including all available features was used. The different 
variations were tested and compared using 8-fold cross validation on the training data with MAPE as the comparison metric. The results show that the features selected through SFS outperform the other 
variations in all cases - see Figure~\ref{fig:inputFeatureComparison}. It is also interesting to note that quantiles generally perform better than mean and standard deviation. For the remainder of the ML study the SFS features are used as input features. 

\section{Results of testing}

\subsection{Regressor: local slightly better than global}

\begin{figure*}
    \begin{subfigure}{0.5\textwidth}
    \includegraphics[width=\textwidth, keepaspectratio]{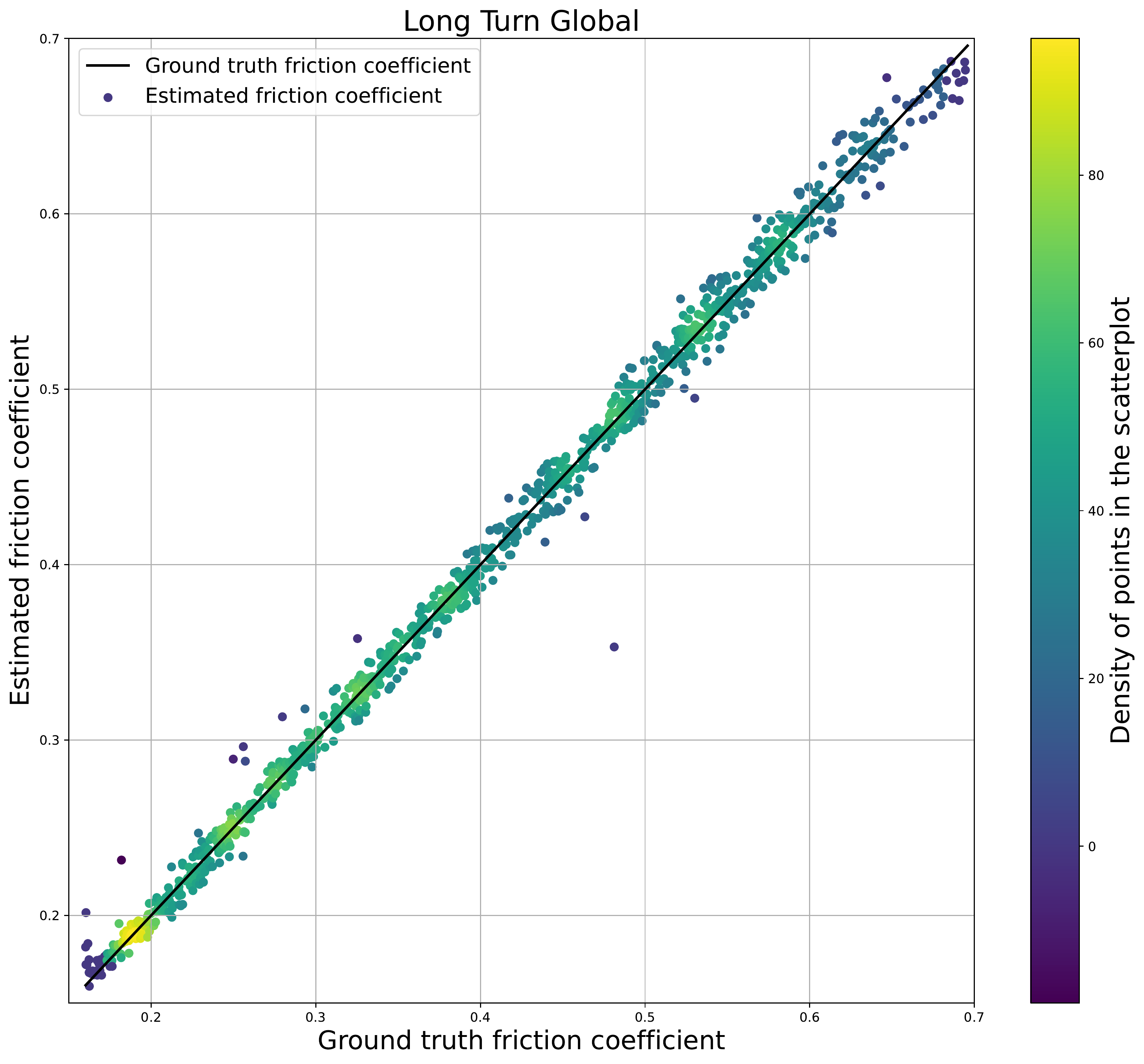}
    \end{subfigure}
    \begin{subfigure}{0.5\textwidth}      \includegraphics[width=\textwidth, keepaspectratio]{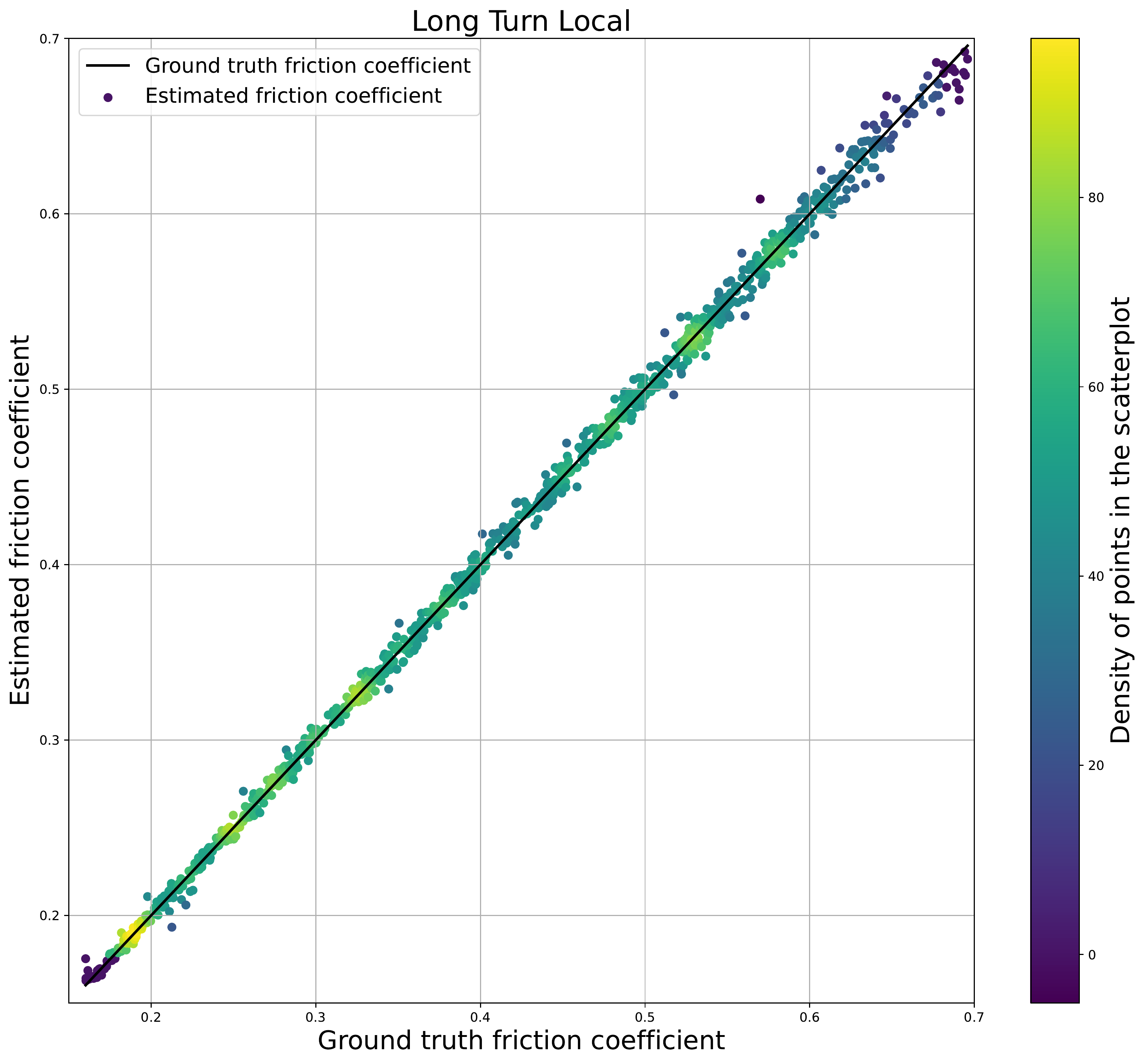}
    \end{subfigure}
    \\
    \begin{subfigure}{0.5\textwidth}
    \includegraphics[width=\textwidth, keepaspectratio]{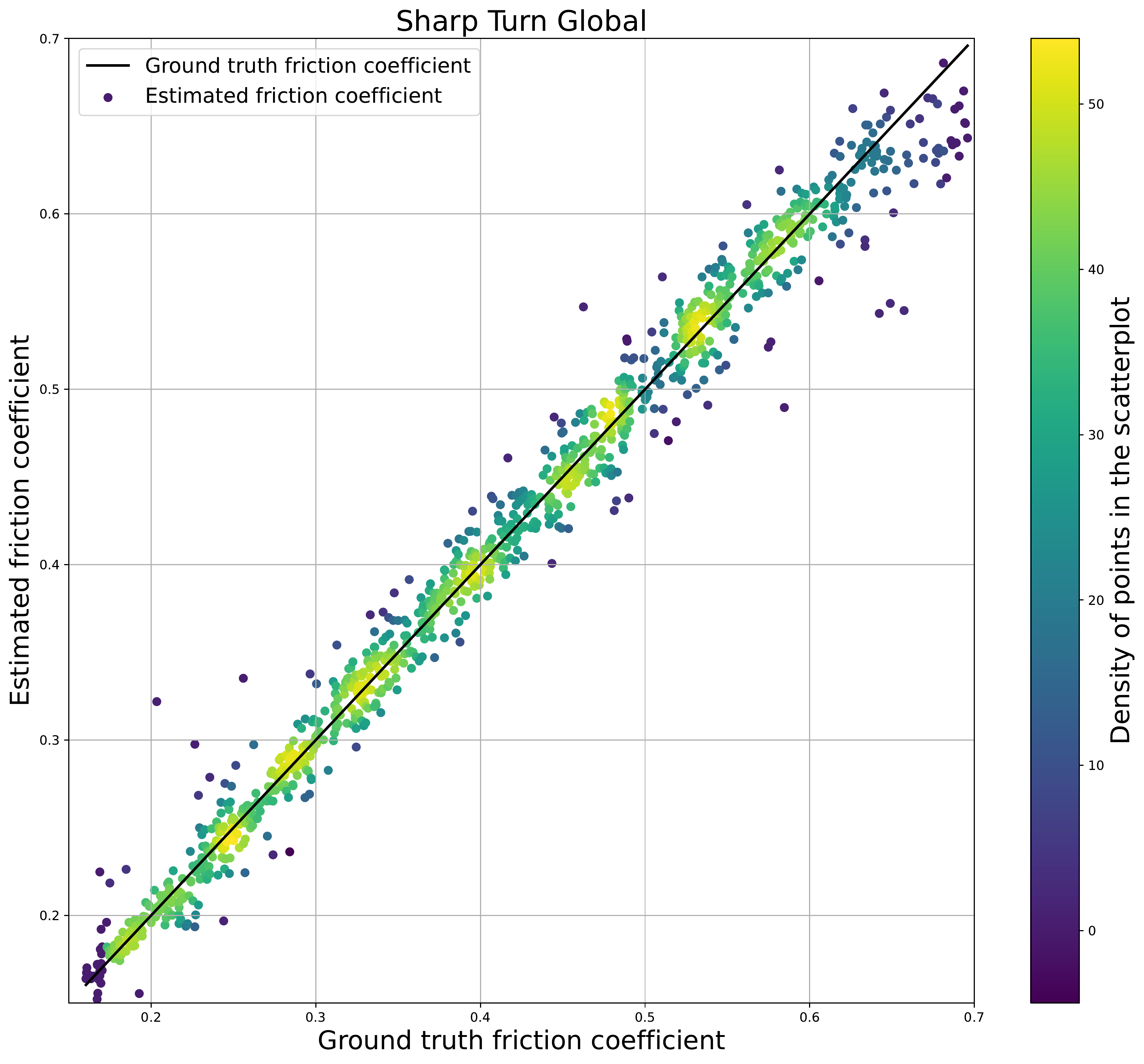}
    \end{subfigure}
    \begin{subfigure}{0.5\textwidth}      \includegraphics[width=\textwidth, keepaspectratio]{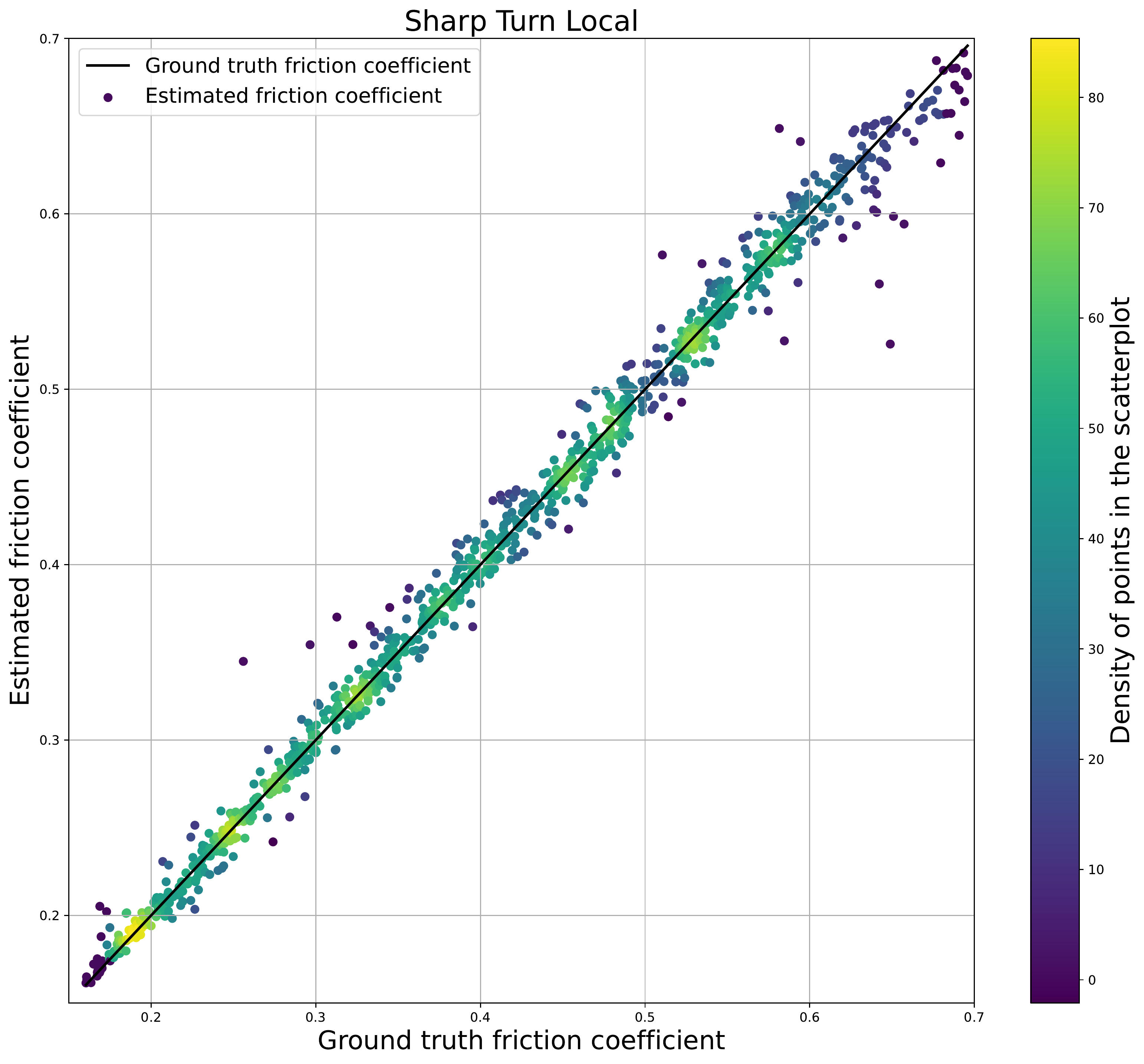}
    \end{subfigure}
\caption{Scatter plots of the estimate versus ground truth friction coefficient, for Long  gentle turn~(top row), and sharp turn~(bottom row). The plots on the left are for the global regressor, and the ones on the right are for the local regressors}
\label{fig:scatterPlotsLongAndSharp}
\end{figure*}

\begin{figure*}    
        \begin{subfigure}{0.5\textwidth}
    \includegraphics[width=\textwidth, keepaspectratio]{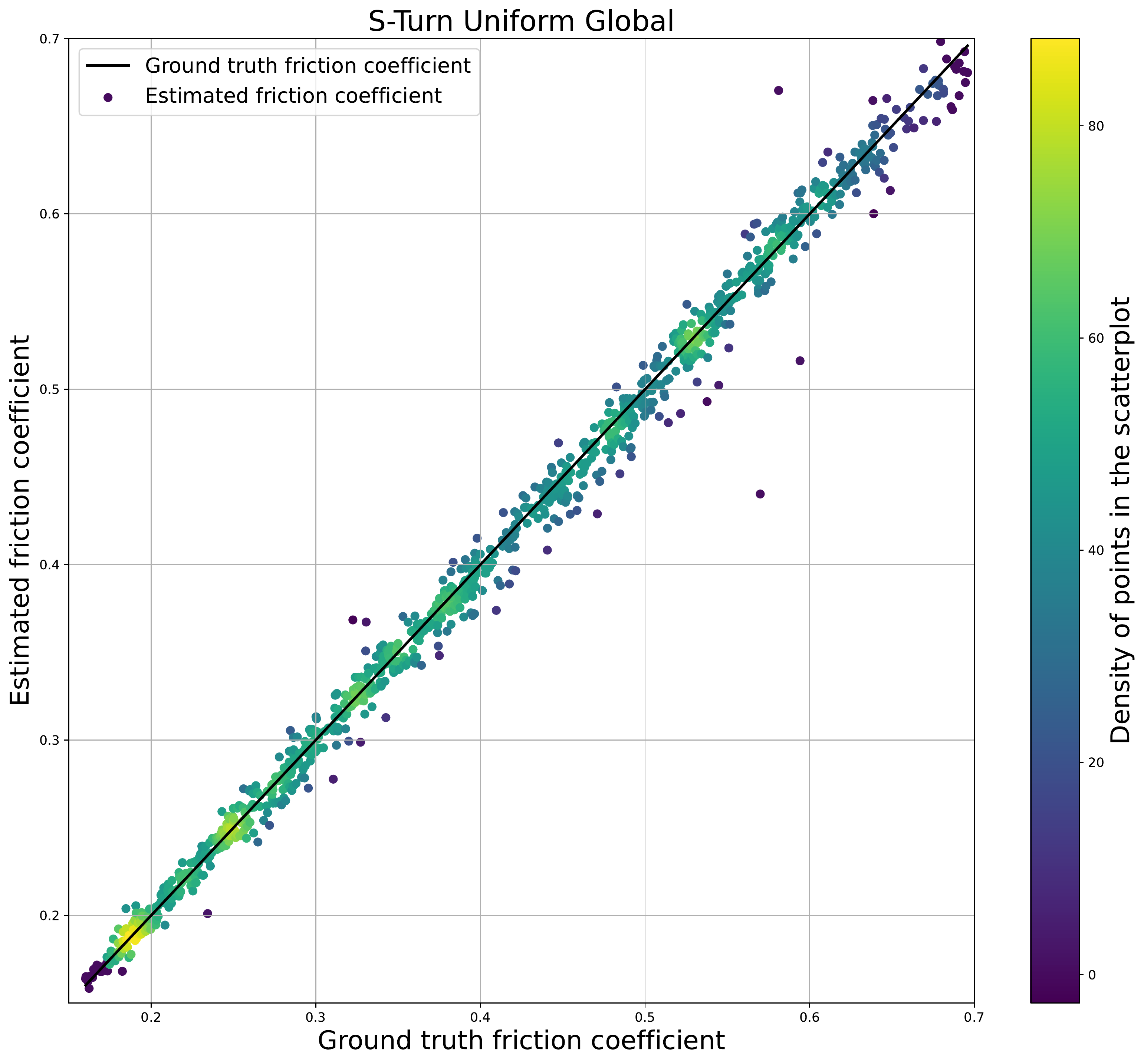}
    \end{subfigure}
    \begin{subfigure}{0.5\textwidth}      \includegraphics[width=\textwidth, keepaspectratio]{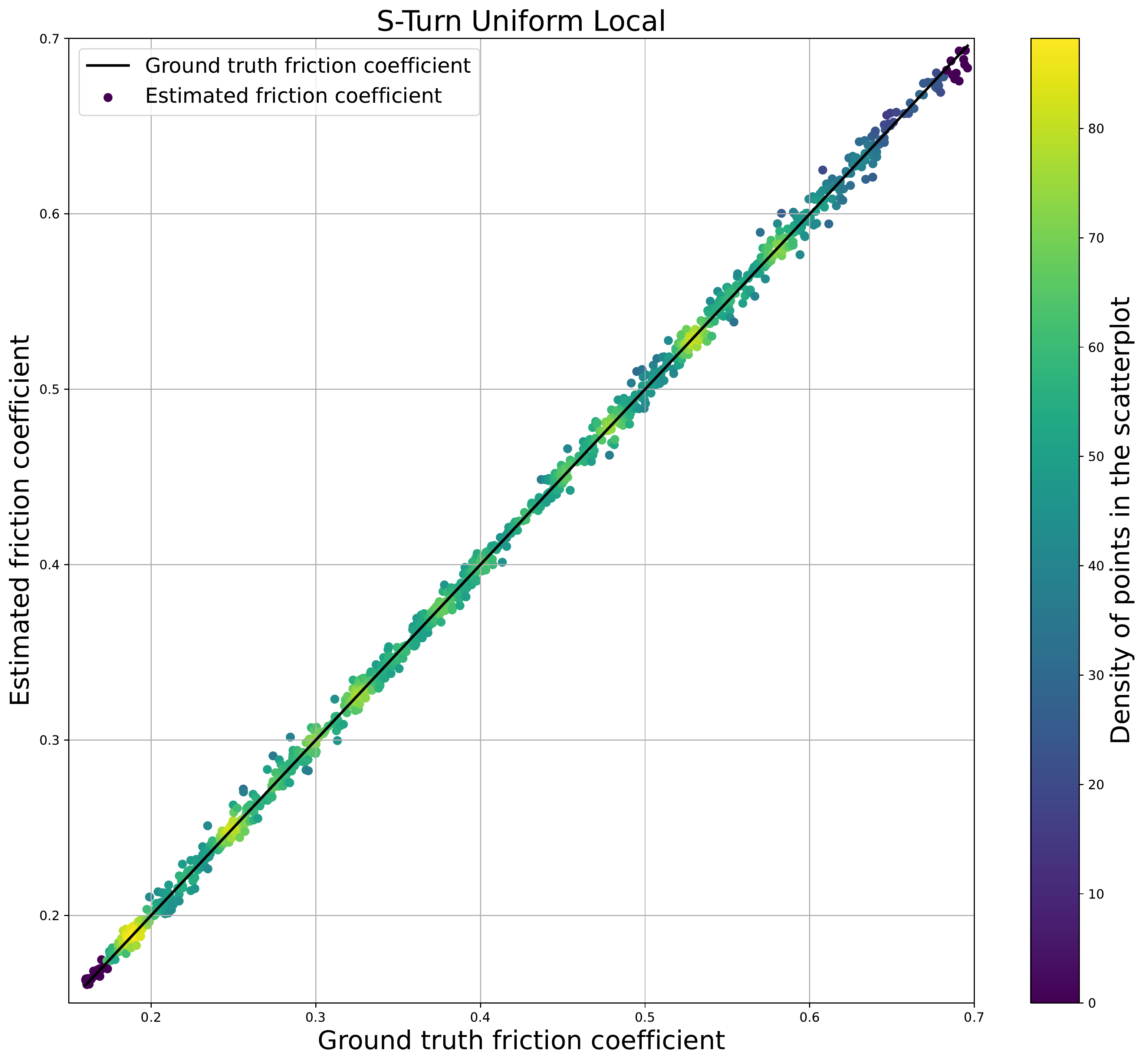}
    \end{subfigure}
    \\
     \begin{subfigure}{0.5\textwidth}
    \includegraphics[width=\textwidth, keepaspectratio]{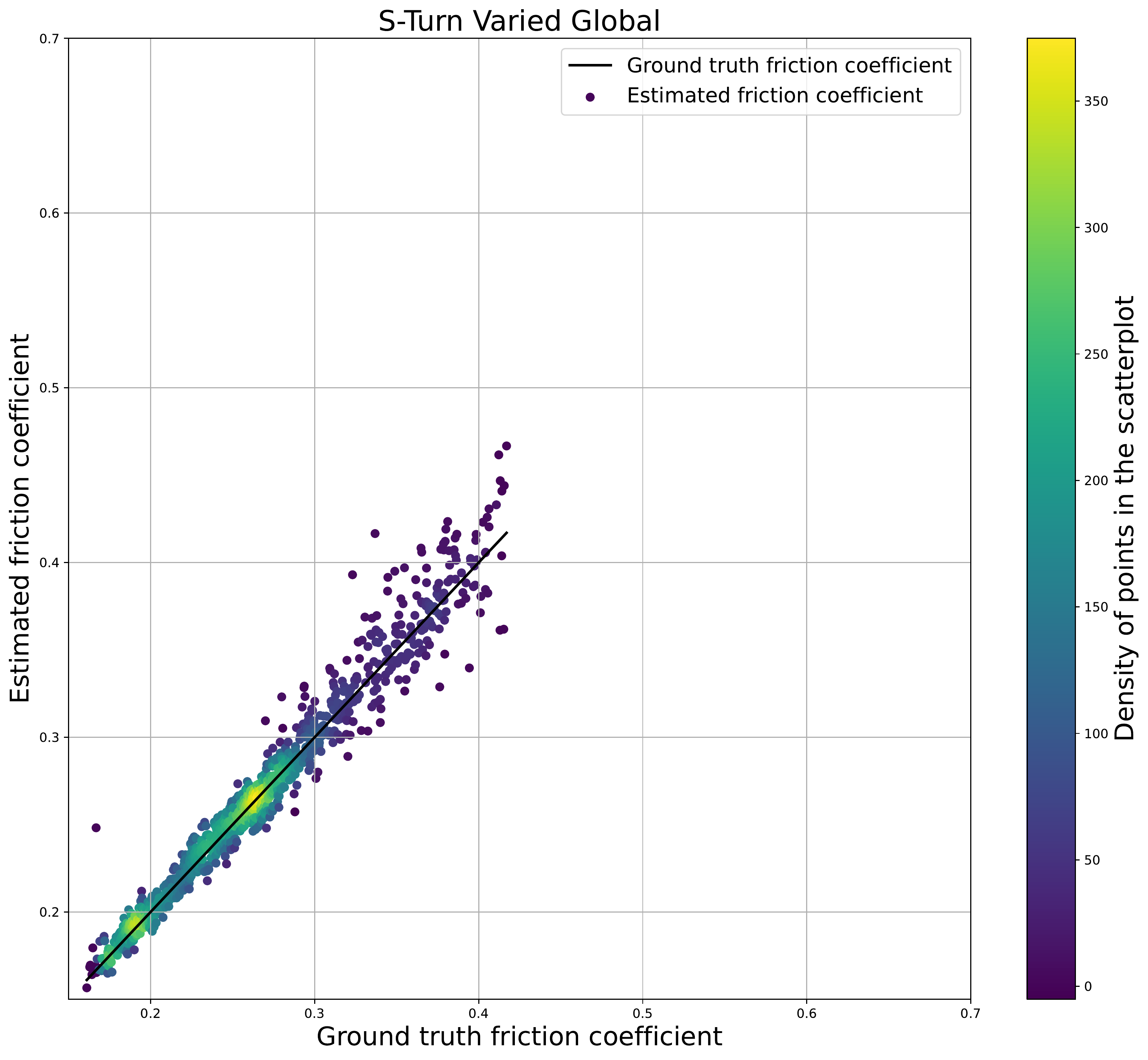}
    \end{subfigure}
    \begin{subfigure}{0.5\textwidth}      \includegraphics[width=\textwidth, keepaspectratio]{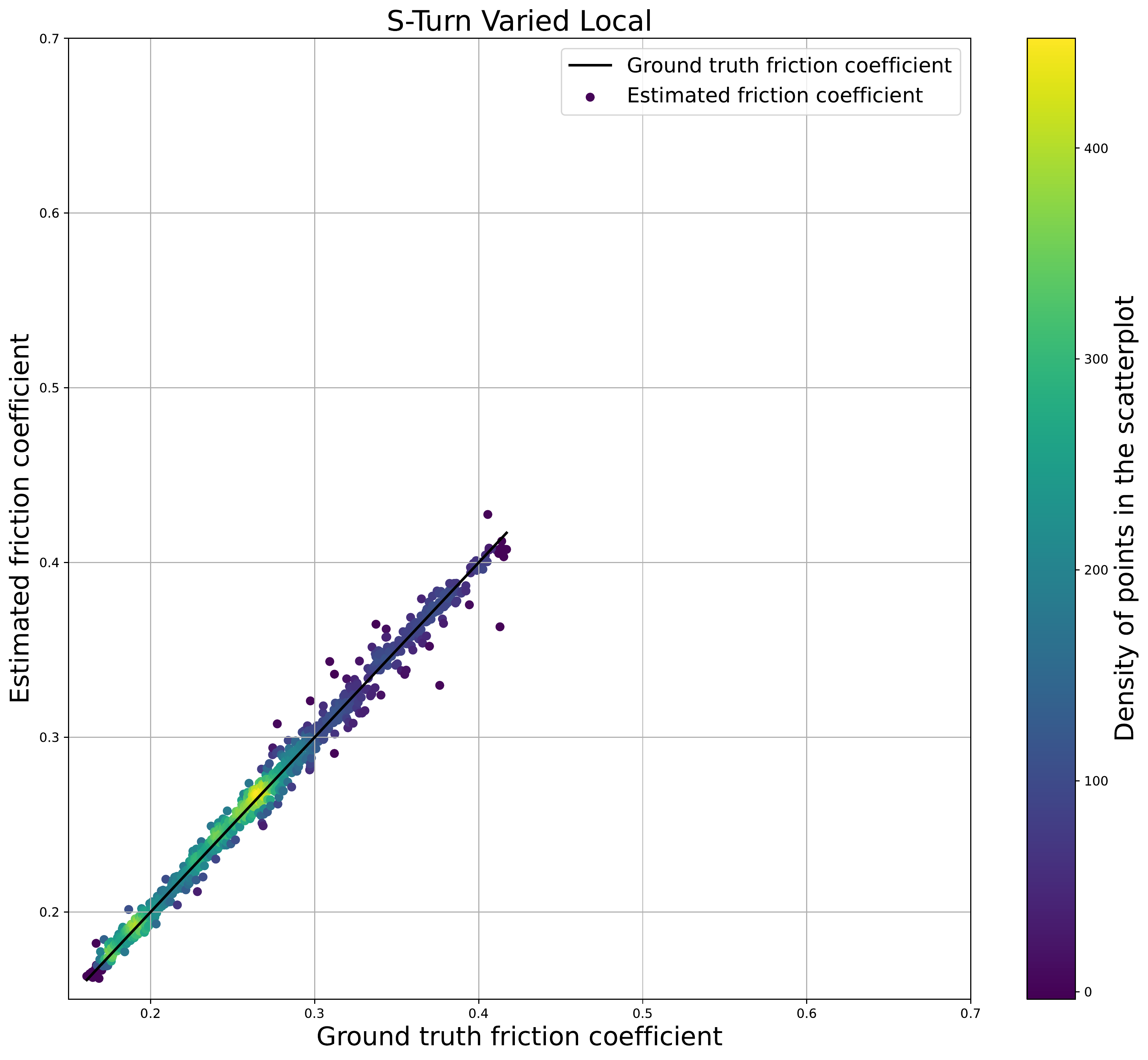}
    \end{subfigure}
    \caption{Scatter plots of the estimate versus ground truth friction coefficient, for S-turn with uniform friction throughout the S~(top row), and two different values of the friction coefficient in the two curved parts of the S~(bottom row). In the bottom row, ground truth is assigned to be the lesser of the two coefficient values. The plots on the left are for the global regressor, and the ones on the right are for the local regressors}
\label{fig:scatterPlotsSturn}
\end{figure*}


\subsubsection{Features from the road geometry}
The global regressor has to be given a means  to  consider the vehicle's kinematic summary in the context of how tightly curved the road section is. 
Hence we need to include features capturing the road geometry.

We add the following road geometry features: length of the road section, the road curvature profile, and a road difficulty index~(RDI) profile. 
The RDI is an index that captures the difficulty caused by insufficient banking at points with high curvature. 
The RDI at each point on  the road segment is calculated as the dimensionless number $ \tan \left(\theta_{\text{ideal}}-\theta_{\text{actual}} \right) , $ 
where $\theta_{\text{actual}}$ is the actual road bank angle, and $ \theta_{\text{ideal}}$ is the uncapped bank angle of Equation~\eqref{eq:bankAngle}. 


The road geometry profile features are statistically summarized in the same manner as the kinematic signals. 

\subsubsection{MAPE performance comparison}
\begin{figure}
{\includegraphics[width=0.475\textwidth, keepaspectratio]{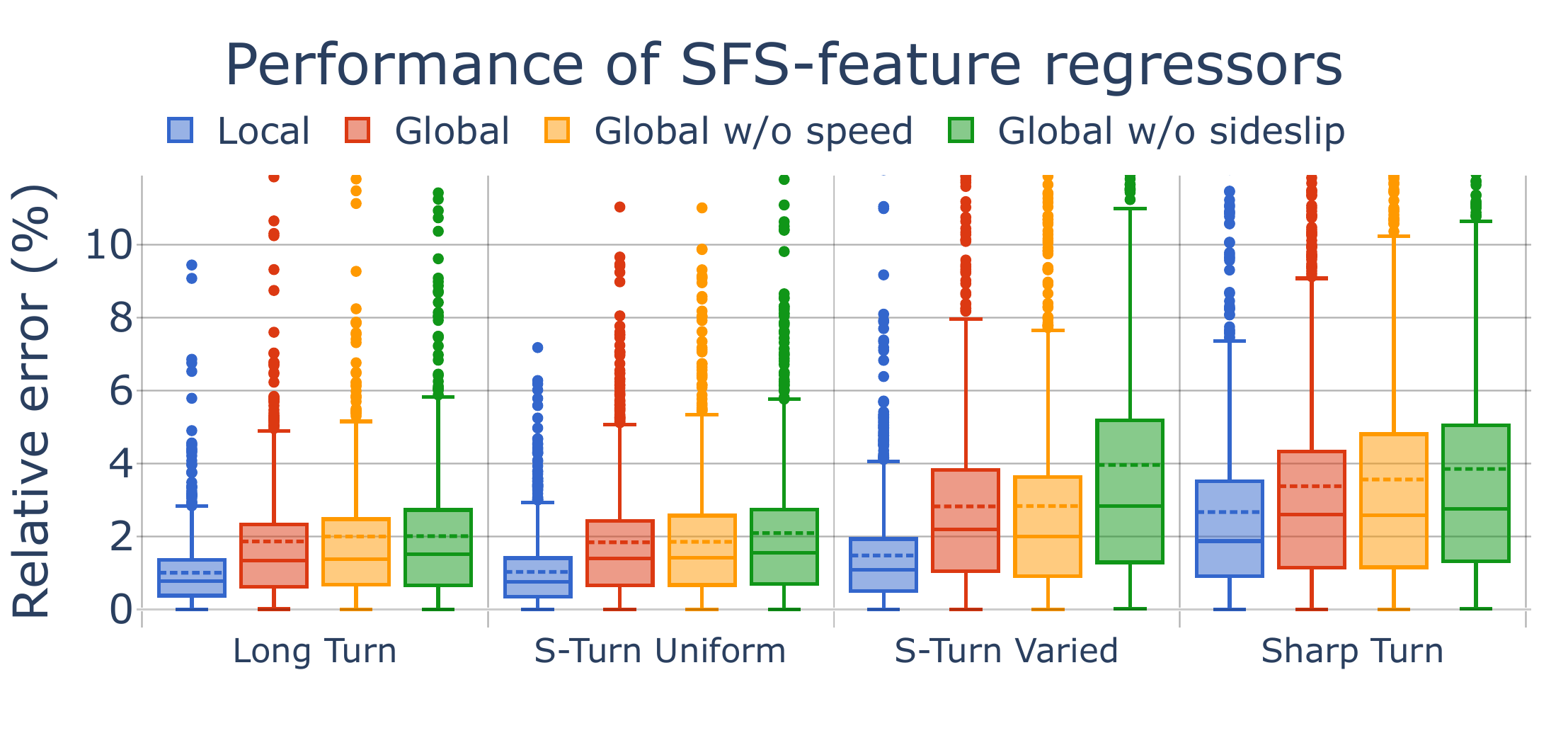}}
\caption{Local and global regressors, on all normal scenarios.}
\label{fig:regressorComparison}
\end{figure}

For each scenario a local XGBoost regressor was fitted to that scenario's training data. And a global XGBoost regressor was fitted to the combined training data of all four scenarios. The performance of the local and global regressors 
were compared, using MAPE as the metric - see Figure~\ref{fig:regressorComparison}. The local regressors perform better than the global regressor, but only marginally so. 
The global regressor's  worst case MAPE  is~3\%, while the local regressors' worst case MAPE is~2\%. Figure \ref{fig:scatterPlotsLongAndSharp} and \ref{fig:scatterPlotsSturn} show the performance of the local and global regressors on each scenario.  

\begin{figure}[!htb]
\begin{center}
    \includegraphics[width=0.475\textwidth, keepaspectratio]{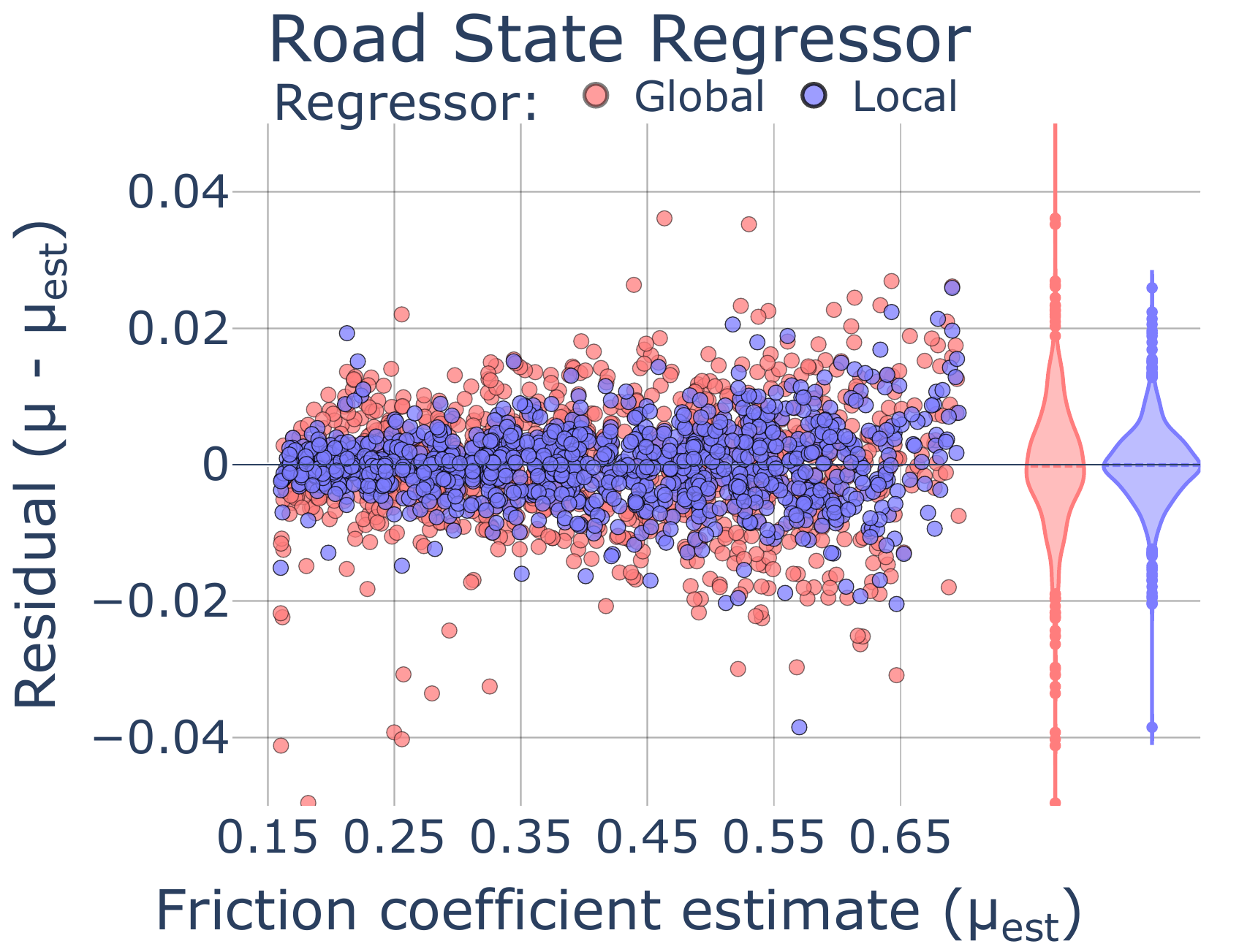}
\end{center}
    \caption{The performance of the local and global regressor on all four scenarios test data.}
  \label{fig:allScenarioResidualPlot}
\end{figure}

We also generated scatter plots of the residual error, for each driving scenario~(see Figures~\ref{fig:scatterPlotsLongAndSharp},~\ref{fig:scatterPlotsSturn}). The plots from the two most complex scenarios can be seen in Figure~\ref{fig:residualPlots} in the appendix. And a combined residual plot for all scenarios can be seen in Figure~\ref{fig:allScenarioResidualPlot}. 
These plots show that
the residuals appear randomly scattered around zero with no clear patterns. 
Additionally, the distribution of the residuals is almost Gaussian, and quite symmetrical. Both are good indicators of 
the soundness of our regressors, for normal driving scenarios.

\subsection{Accuracy of the RSU's interval estimate}
The interval of friction coefficients at a road section captures the variation in the friction coefficient, which individual
vehicles experience due to differences in tyre wear, tyre pressure, load variation etc.
As described in Section~\ref{section:intervalsOfFrictionConditions}, the friction coefficients of  individual vehicles are assumed to be symmetrical about the midpoint of the interval of friction that is in operation. Then
the task of estimating the interval reduces to estimating its midpoint. As declared before, we set the RSU's estimate of this interval midpoint to be the median of individual
friction coefficient estimates from the
batch of vehicles that have recently passed.



\begin{figure}
{\includegraphics[width=0.40\textwidth]{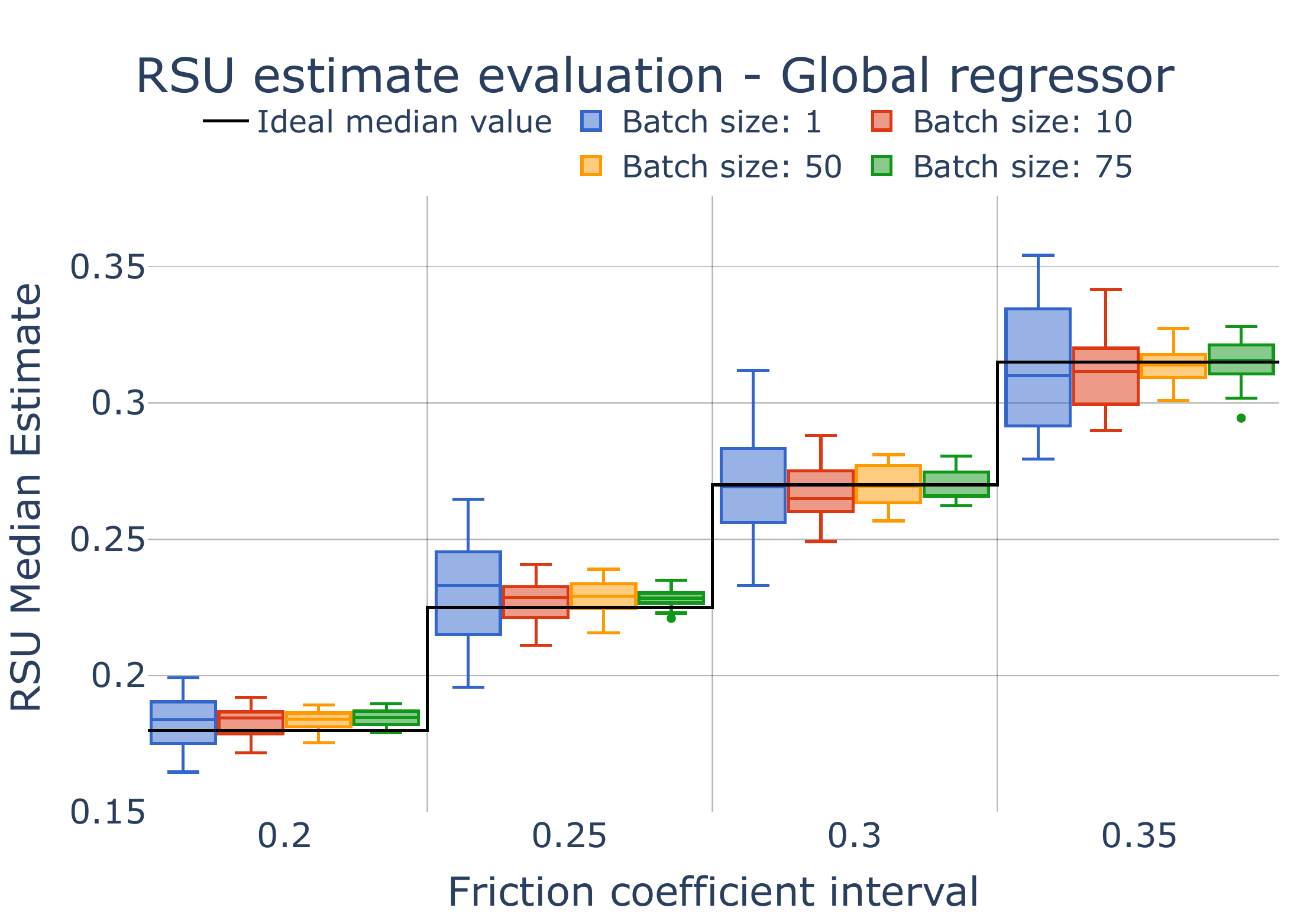}}
\caption{Accuracy of the RSU's interval estimate in normal scenarios, as a function of traffic intensity at the road section.}
\label{fig:rsuEstimateBatchSizes}
\end{figure}

Figure \ref{fig:rsuEstimateBatchSizes} shows how the  batch size affects the performance of the collaborative regressor. Using a batch size of~1 is the same as using a single estimate.
As the batch size increases, the RSU's median estimate approaches the true midpoint of the friction coefficient interval.
We get a barely tolerable worst case error with a batch size of 10~vehicles, and with a batch size of at least 50~vehicles we get an acceptable worst case error of about 10\%, which occurs in the band with the base friction coefficient value at~0.35.

\subsection{Comparison with a commercial service}
A bachelor's thesis~\cite{andreasson2017luleaBachelorThesisNiraDynamics} documents the only  information in the public domain on the accuracy of the  friction coefficient estimation system provided
by the company Nira Dynamics~\cite{niraDynamics}.
The road test records the Nira dynamics system's estimates, together with the ground truth friction coefficient values measured by a custom-made  friction wheel. The tests were run by driving a single vehicle on snowy and icy road sections, each of which offered variable friction coefficients.

Twenty-one test runs were reported in Table~6.1 of the thesis. For each run, we get: (a)~the RMS value of the residual error, and (b)~the average value of ground truth friction coefficient.
From these values, we can estimate that relative to the ground truth, the RMS error is on average 14.81\%, and in the median 11.33\%.

Our system is at least twice as good as the above-mentioned commercial system.  For our system, the RMS error figures for the regressors  estimating the friction coefficient for individual vehicles are no more than 5\% of the ground truth. Furthermore, in our system, if the RSU's interval estimate is based on at least fifty vehicles, then we get a worst case error that lies between 6\% and~10\%.

\subsection{Testing with extreme manoeuvres, unseen in training}
So far, we have used normal driving scenarios at curved road sections, to generate training data for our regressor, and also to generate the test data that was used in the previous section. Those driving scenarios are normal in that sense that: (i)~the considered vehicle speeds are within about twenty percent of the rated speeds, and (ii)~the road geometries do not demand extreme manoeuvres - for example, in our training dataset the sideslip signal is less than 3~degrees.

Extreme manoeuvres have been considered by some previous works~\cite{lampeZiaukas2022lstmRoadFriction,hanyang2018annFriction}. These works have considered manoeuvres where the sideslip and steering angles are an order of magnitude higher than in our training dataset.
For such extreme manoeuvres in practice, it is hard for the road authority to generate training data. For such manoeuvres, it is not possible to use the specialized vehicles~\cite{viaFriction} for measuring road friction, because of unpredictability in the the motion and functioning of the measurement trailers.

\begin{figure*}
    \begin{subfigure}{0.49\textwidth}
        \includegraphics[width=\textwidth, keepaspectratio]{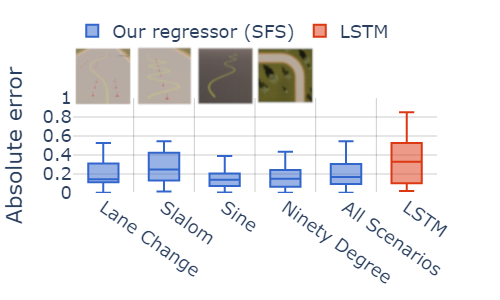}
    \end{subfigure}
    \begin{subfigure}{0.49\textwidth}
        \includegraphics[width=\textwidth, keepaspectratio]{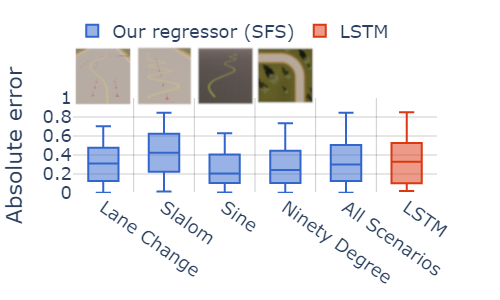}
    \end{subfigure}
\caption{Box plots of the error in friction coefficient, for test on the unseen data from extreme manoeuvres. Left: friction ranging from 0.2 to 0.7. Right: friction ranging from 0.2 to 1. Note: the box plot of the LSTM~\cite{lampeZiaukas2022lstmRoadFriction} is for the RMS errors and 0.2 to 1 friction range.}
\label{fig:boxPlotsUnseenData}
\end{figure*}

Hence we test our regressors on data from extreme manoeuvres, but without having seen any such data in our training dataset.
We generated test data from the extreme manoeuvres  described in~\cite{lampeZiaukas2022lstmRoadFriction}: slalom, tight lane change, ninety degree turn at relatively high speed, sine wave driving. The friction coefficient ranges from~0.2 all the way up to~1.

Firstly, the performance of our regressor becomes an order of magnitude worse, than its performance
on the normal driving scenarios at curved sections.
In particular, our regressor vastly underestimates the friction coefficient, if the ground truth value is above 0.5~(see Figure~11 in~\cite{langstrandRabi2023regressorForFrictionCoefficient}).
But our regressor has a slightly smaller RMS error than the regressor of ~\cite{lampeZiaukas2022lstmRoadFriction,hanyang2018annFriction} (compare Figure~\ref{fig:boxPlotsUnseenData} with Figure~6(b) of~\cite{lampeZiaukas2022lstmRoadFriction}). This is only a superficial comparison, because Figure~6(b)~\cite{lampeZiaukas2022lstmRoadFriction}  gives the box plots of RMS errors, wheras we use the absolute error.

\begin{figure*}
    \begin{subfigure}{0.49\textwidth}
        \includegraphics[width=\textwidth, keepaspectratio]{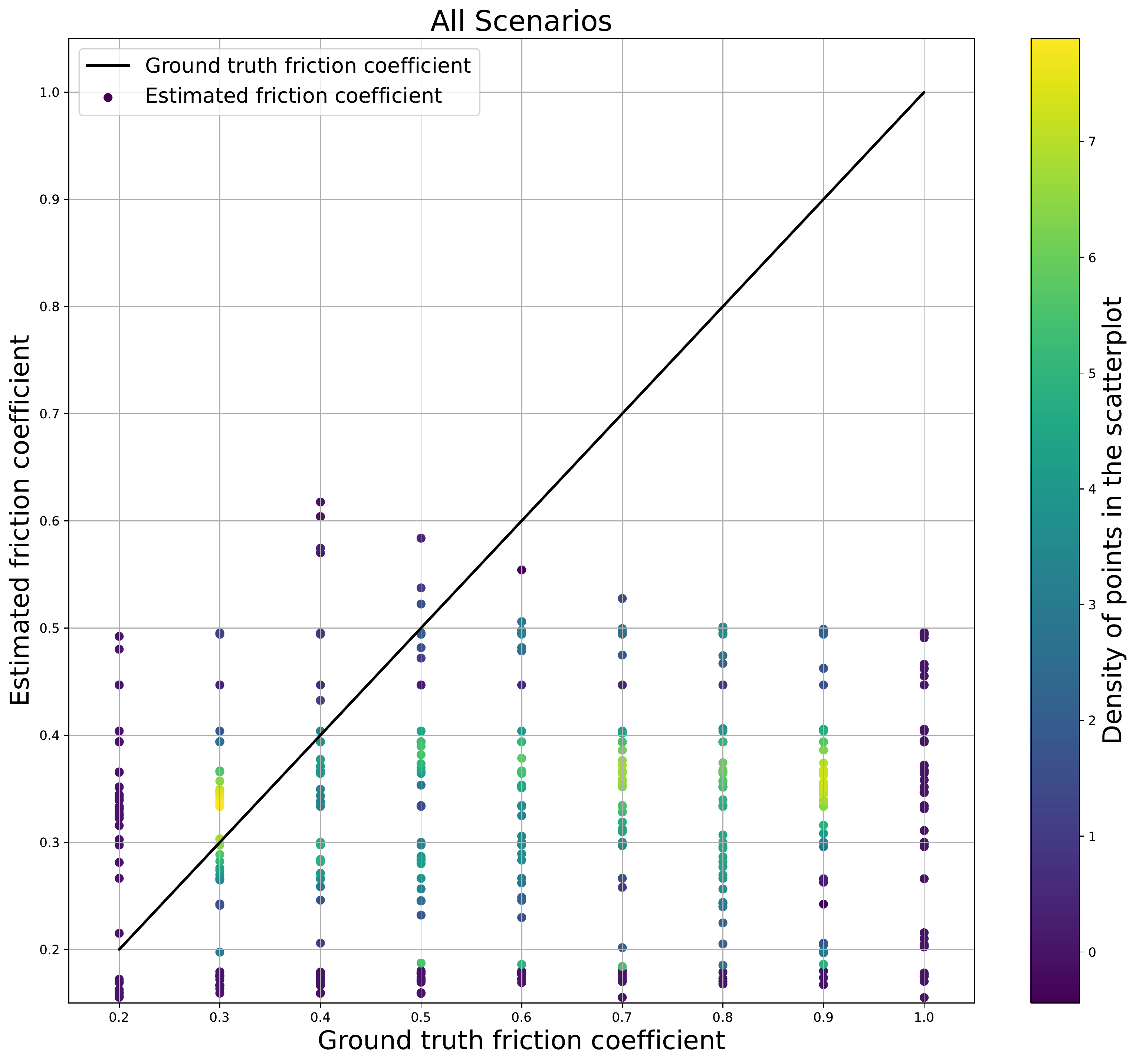}
    \end{subfigure}
    \begin{subfigure}{0.49\textwidth}
        \includegraphics[width=\textwidth, keepaspectratio]{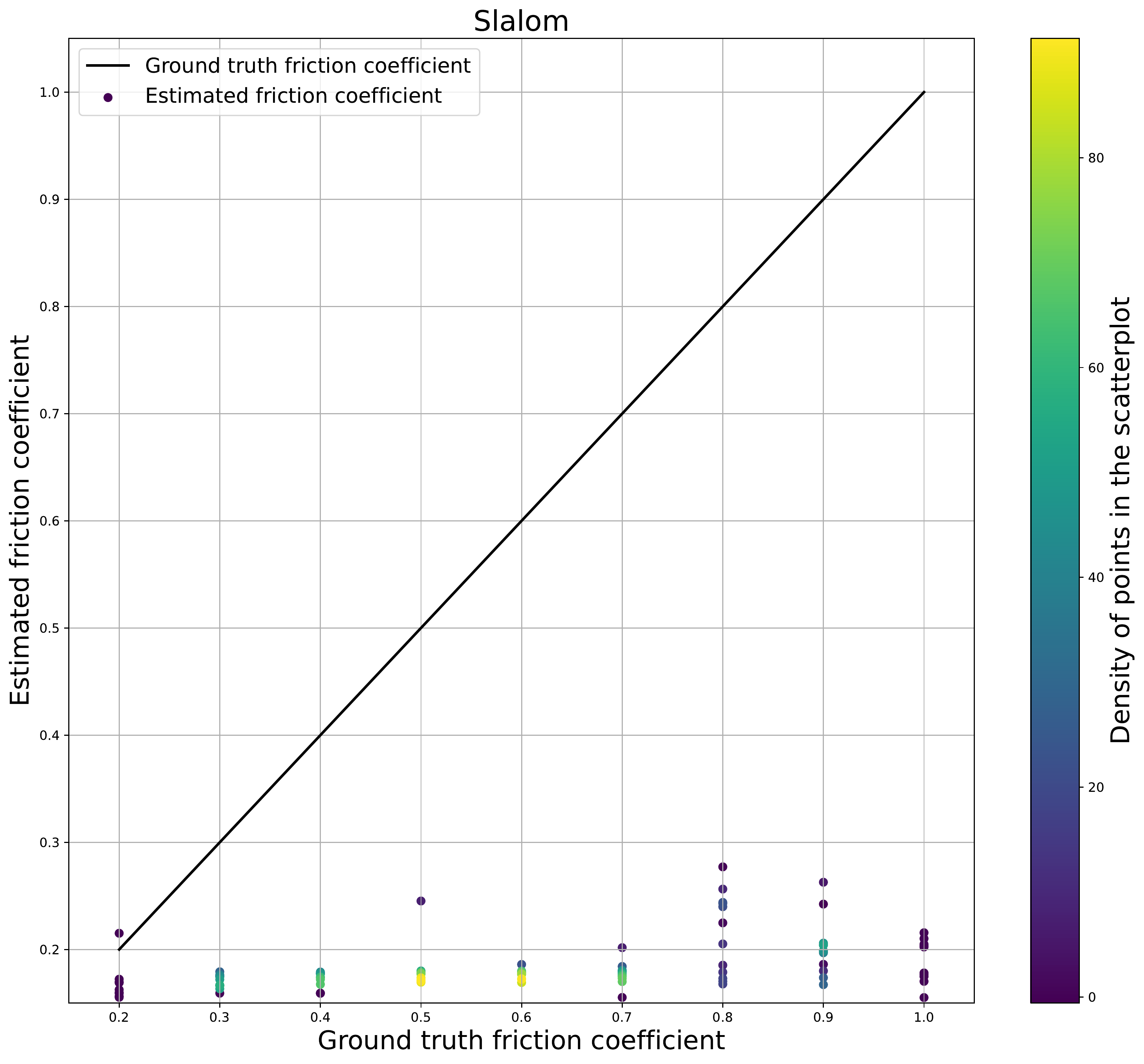}
    \end{subfigure}
    \caption{Scatter plots of the estimate versus ground truth friction coefficient, for all the extreme manoeuvres~(left), and the slalom manoeuvre in particular~(right). The global regressor is used, as there are no locally trained regressors for these extreme road sections.}
\label{fig:linePlotsUnseenData}
\end{figure*}

\begin{figure}
    \includegraphics[width=0.48\textwidth, keepaspectratio]{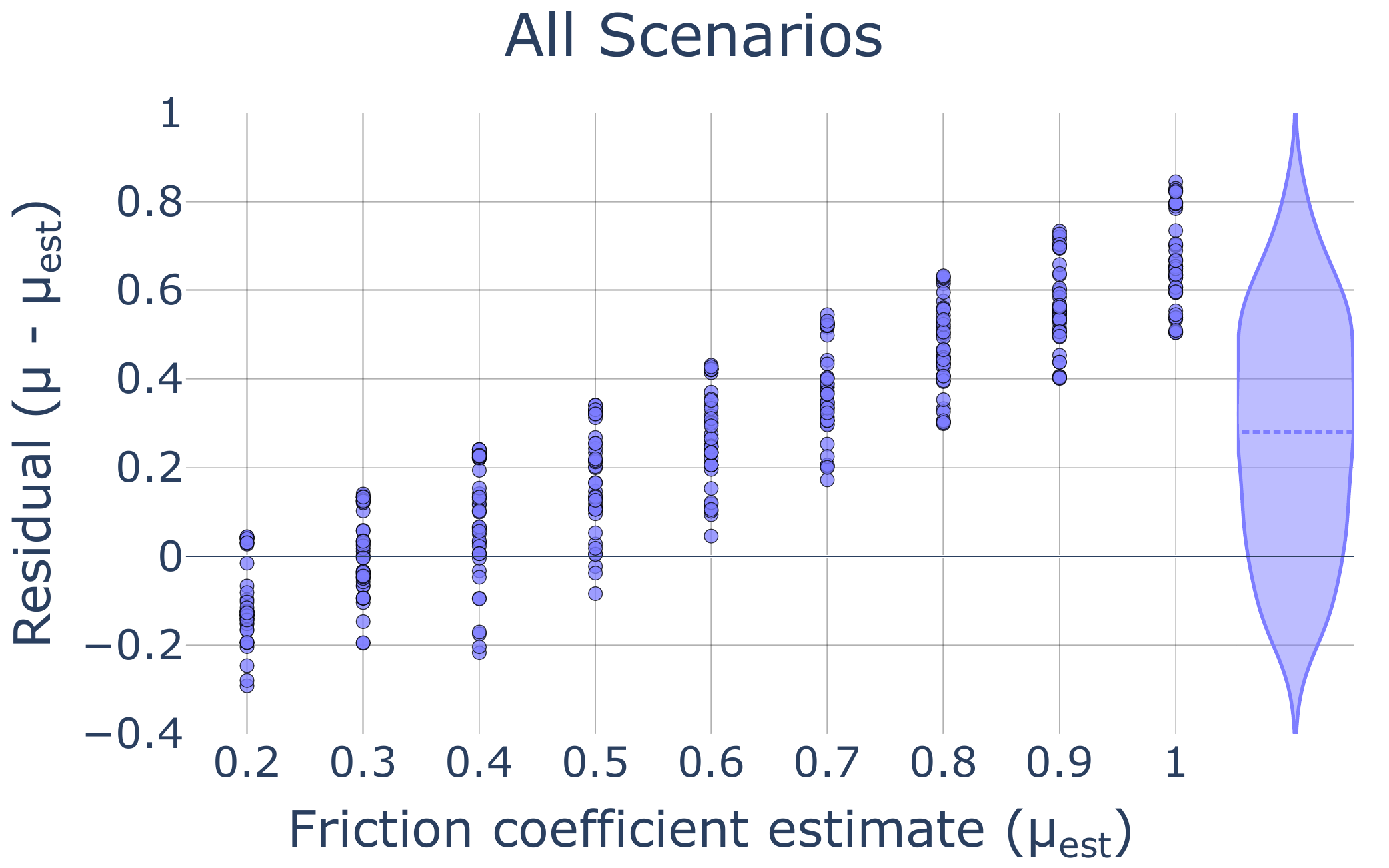}
    \caption{The residual plot from our regressor's (SFS) friction estimate on all the unseen, extreme manoeuvres.}
\label{fig:residualPlotsUnseenData}
\end{figure}

\section{Summary}

An XGBoost regressor with Random forests as the underlying tree algorithm 
performed the best on our datasets. 
We tried different input feature variations and found that a subset of the available features selected by a SFS algorithm had the best performance. Out of these features the steering angle, sideslip angle, yaw rate excess and the vehicle speed were the most important. This drastically reduces the amount of input features from 82 to 15 and therefore the complexity of the resulting ML model as well. 

The local regressors performed only marginally better than the global regressor. This news is good, because there are practical benefits to using a global regressor. These include simpler deployment, maintenance and continuous improvement of the regressor as new data is collected. Crucially, the global regressor will be trained on a much larger and varied dataset when compared to a local regressor, which should result in better generalisation to unseen data.

Surprisingly, the road geometry features rank low and barely make it into the list of selected features.

Perhaps unsurprisingly, our regressor performed  quite badly, when tested  with data from extreme driving manoeuvres
that were unseen in the training data.



The RSU's median estimate of the friction coefficient robustly approaches the true midpoint of the friction interval, as the batch size of vehicle summaries increases. 
Since broadcasting friction coefficient estimates impinges on safety, we recommend using RSU interval estimates only if we have batch sizes of at least 50~vehicles.

\subsubsection{Can we ignore the sideslip?}
The sideslip is the most challenging vehicle state to estimate. If we were to remove it from our kinematic summaries, then there shall be a large drop in the complexity of designing and executing the onboard processing 

Figure~\ref{fig:regressorComparison} shows
the effect on the accuracy of individual friction coefficient estimates, if we drop features connected to the sideslip and its rate - see the box plots in green, in comparison to the others.


\subsubsection{Shortcut learning has been avoided}
There was a risk of shortcut learning~\cite{tubingen2020shortcutLearning},
where regressors could take cues from the speed profile rather than the
`true' kinematic signals, namely the sideslip and the yaw rate excess.
That is to say that regressors learn to correlate changes in speed caused by the vehicle driver model on the one hand, with the coefficient of friction on the other. We tested for this by eliminating speed related features, and then comparing the performance to that of the original regressor. The regression performance barely changed - see the orange and red coloured box plots in Figure~\ref{fig:regressorComparison}. Hence there is no indication of shortcut learning. 

\subsection{Gaps in evidence}
Our study relies upon simulated vehicle data. 
We need validation via real vehicle tests. 

Our simulations had only one type of vehicle, with one driver model, and a vehicle state observer tuned to that vehicle. Hence one may conclude that any positive result we show may be limited to that sort of vehicle. If such a pessimistic conclusion turns out to be true, then this shall limit the accuracy of any practical deployment. In specific, if only a single sort of vehicle can be enrolled, then the vehicle batch sizes for collaborative estimation shall dramatically decrease. Then instead of the road authority running the friction estimation service, it may be easier for the relevant vehicle manufacturer to run the virtual RSUs.




\section*{Acknowledgement}
This work was supported by the Norwegian research council through the project, {\textit{CriSp: Finding a CRitical SPeed function ahead of a road section for vehicles in motion,}} with the project number~302327, and was also supported by {\O}stfold University College through its funding for research time and purchase of the Dyna4 simulator. We are grateful for this support. We also thank Jim Tørresen and Sabita~Maharjan for valuable feedback on an earlier draft of this manuscript.
Figure~\ref{fig:theCriSpStack} is an imitation of Figure~1 from~\cite{autowareOnBoard2018iccps}.
Figure~\ref{fig:vehicleRSUMessaging} uses background images from \url{Depositphotos.com}.

\bibliographystyle{acm}
\bibliography{mabenRefs,jpRefs}

\section*{APPENDIX}
Here we collect segments of explanations, configurations and graphs that support the descriptions and conclusions in the main text.
\subsection{Options for Wireless communications%
\label{section:wirelessStandardDesignChoice}}
The choice of the V2X communication standard does not affect the behaviour or performance of the system.  Rather it affects deployability in certain geographic areas. There are two factors to consider: the reach of the wireless communication, and the cost of installing and servicing any road side infrastructure that may be needed.

If we were to prescribe WiFi communications to RSUs, then we can deploy the system even at remote locations where there is no cellular coverage. But this requires us to physically install a RSU at every road section where we want our system. With this choice, communications have to happen within time windows that are approximately a minute long. In specific the packet exchanges between a vehicle an RSU have to happen within time windows that occur when the vehicle approaches the road section and is within range of the RSU, and also when the vehicle leaves the road section and is still within range of the RSU.

If we were to prescribe cellular communications, then we can only deploy the system at locations with cellular reception. But at these spots we can make use of the already existing cellular infrastructure (4G or 5G). The communications do not have to happen within narrow time windows - packets can be exchanged some minutes to tens of minutes before and after the vehicle has crossed an identified road section.

\subsection{Extra details on vehicle state estimation}
The ideal, uncapped road bank angle is given by: 
\begin{equation}\label{eq:bankAngle}
    \theta_{\text{ideal}} = \frac{V^2}{127R} - 0.14 
\end{equation}
where  $V$ is the speed limit of the road section, $R$ is the radius of the turn, and $0.14$ is the side friction coefficient.

\subsubsection{Observer tuning parameters}

\begin{singlespace}
\begin{table}[H]
\centering
\small
\caption{
\bf{Tuned observer parameter values}}
\tabcolsep=0.11cm
$\begin{array}{l l l l}
\hline
Parameter & Value & Parameter & Value \\
\hline
\alpha_0 & 15 & \alpha_1 & 14\\
\alpha_2 & 5& A_x^{threshold} & 1\\
\sigma_{d\omega_1} & 4.5 & \sigma_{d\omega_2} & 1\\
\sigma_{\omega_z} & 0.22 & \sigma_{d\omega_z} & 0.22\\
\sigma_\delta & 0.12 & \sigma_{d\delta} & 0.08\\
\sigma_{\dot{\beta}} & 0.05 & \sigma_{d\dot{\beta}} & 0.3
\end{array}$
\label{table:tunedObserverConstants}
\end{table}
\end{singlespace}

\subsection{Extra details on the simulation%
\label{section:simulationExtraDeatils}}

\subsubsection{Vegkart scenario queries \label{section:roadQueries}}
\begin{itemize}[leftmargin=*]
\item[] Long turn: \url{https://tinyurl.com/LongTurnScenarioQuery} 
\item[] S-turn: \url{https://tinyurl.com/STurnScenarioQuery} 
\item[] Sharp turn: \url{https://tinyurl.com/SharpTurnScenarioQuery}
\end{itemize}

\begin{figure}
\centering
\includegraphics[height=.15\textheight, keepaspectratio]{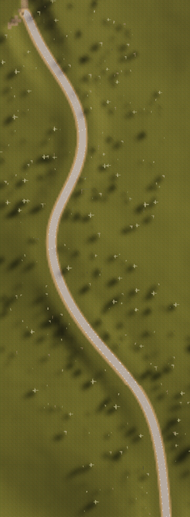}
\ \ \ \ \ \ \ \ \ \ \ \ 
\includegraphics[height=.15\textheight, keepaspectratio]{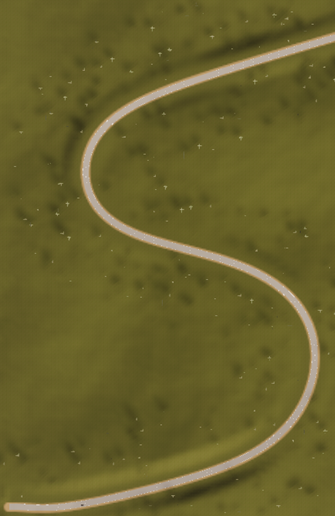}
\ \ \ \ \ \ \ \ \ \ \ \ 
\includegraphics[height=.15\textheight, keepaspectratio]{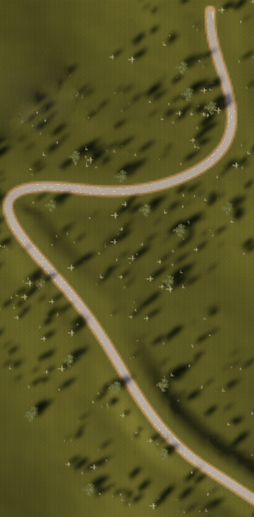}
\caption{Top views of the three normal driving  scenarios. From left to right: Long turn, S-turn, and Sharp turn.}
\label{fig:scenarios}
\end{figure}
\subsection{Extra details on ML performance}

\begin{figure}[H]
\includegraphics[width=0.44\textwidth, height=0.25\textheight, keepaspectratio]{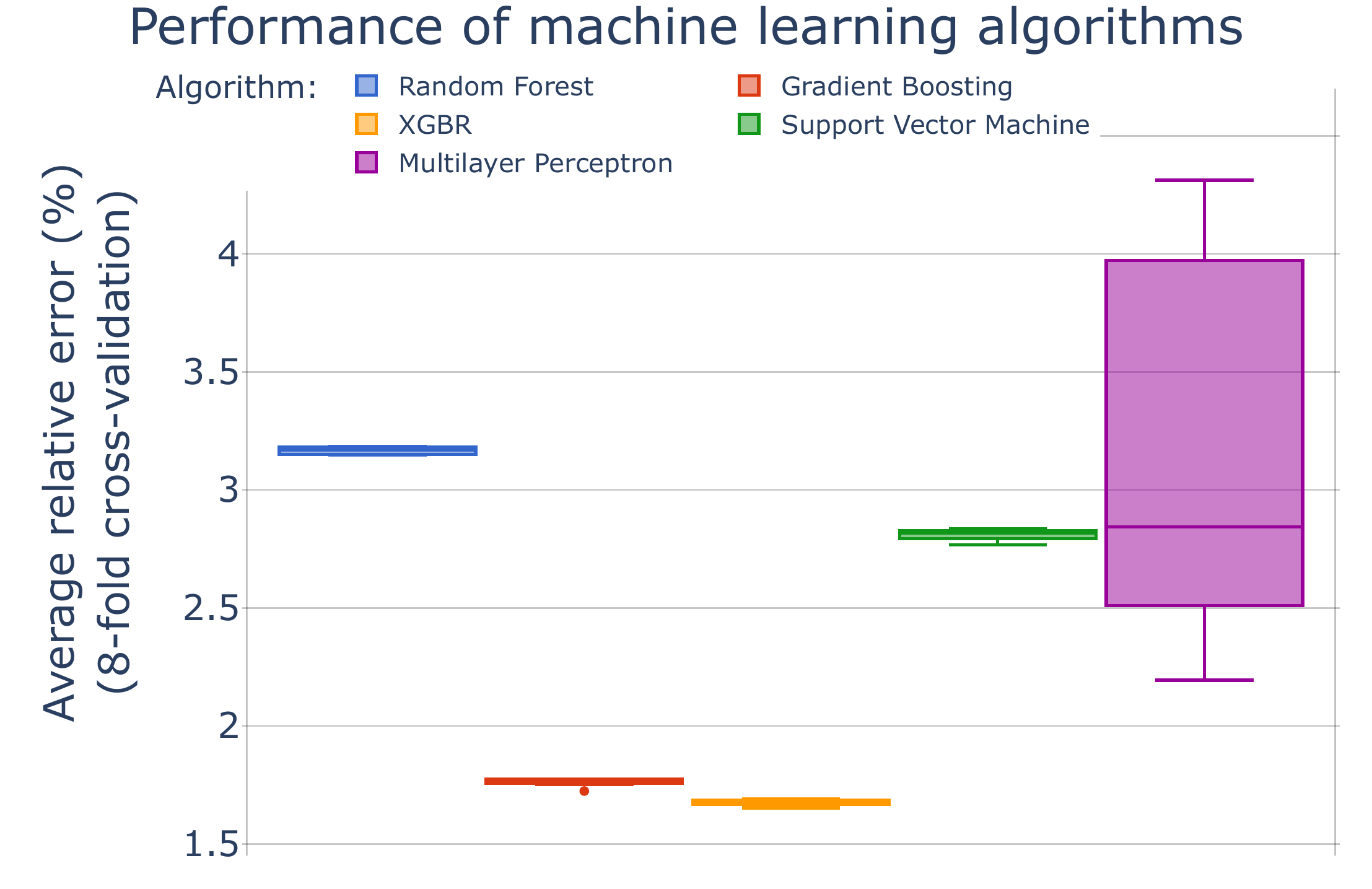}
\caption{The performance of different machine learning algorithms on a combined  training dataset.}
\label{fig:modelCVPerformance}
\end{figure}

\begin{figure}[H]
{\includegraphics[width=0.40\textwidth, keepaspectratio]{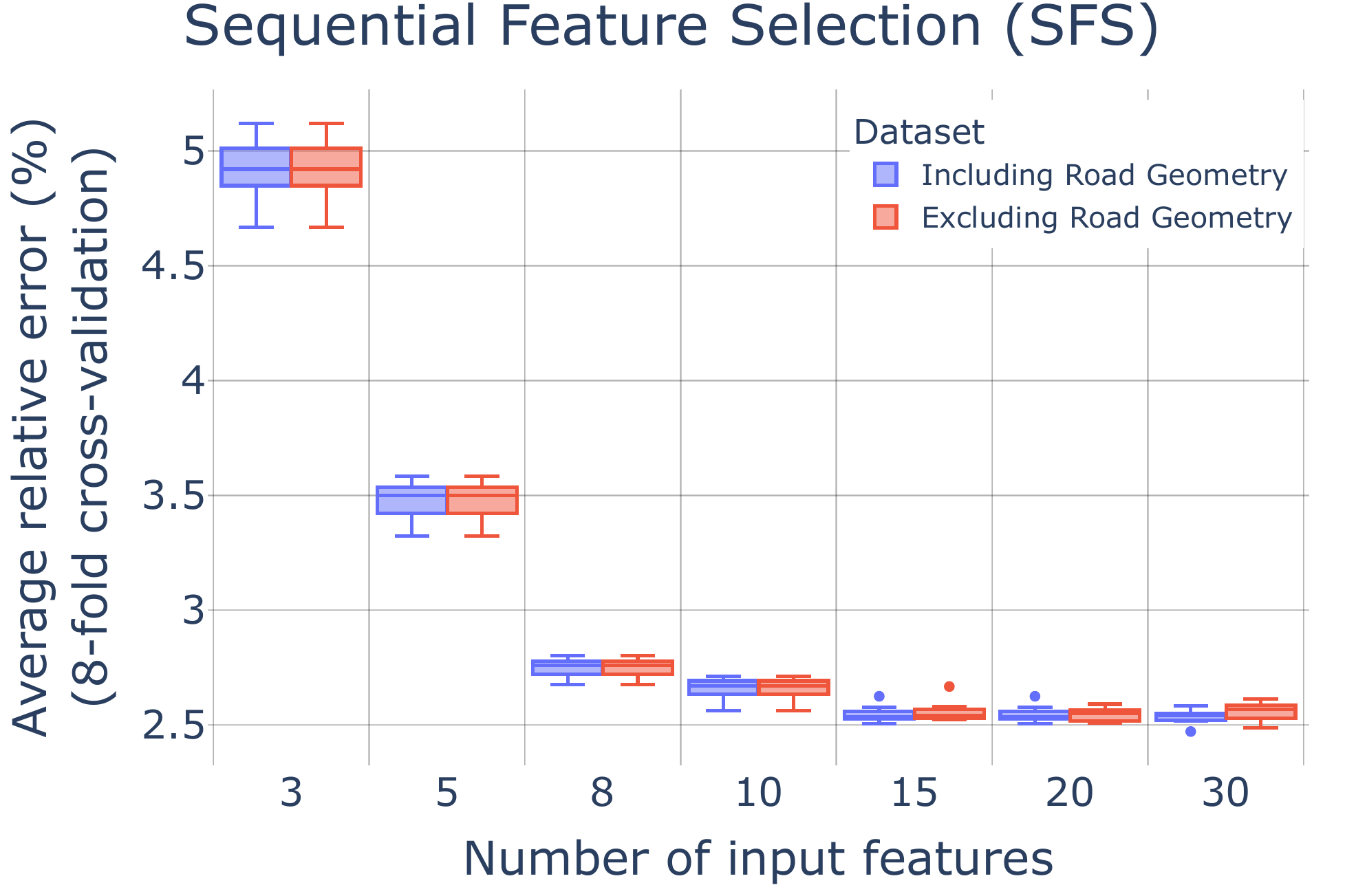}}
\caption{Performance comparison of XGBoost regressors with increasing number of input features.}
\label{fig:nmbFeaturePerformance}
\end{figure}

\begin{figure*}
    \begin{subfigure}{0.43\textwidth}
        \includegraphics[width=\textwidth, keepaspectratio]{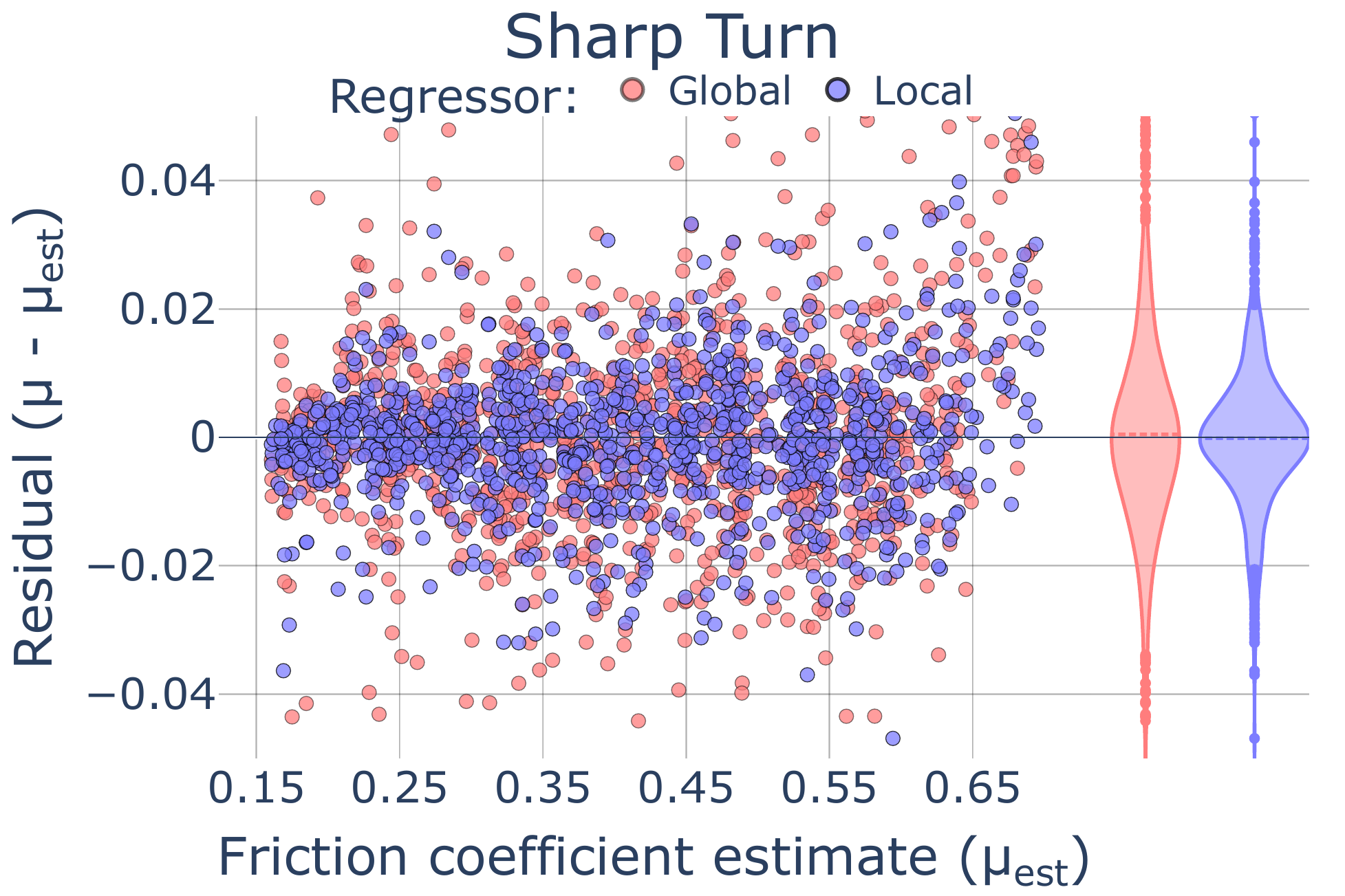}
    \end{subfigure}
    \begin{subfigure}{0.43\textwidth}
        \includegraphics[width=\textwidth, keepaspectratio]{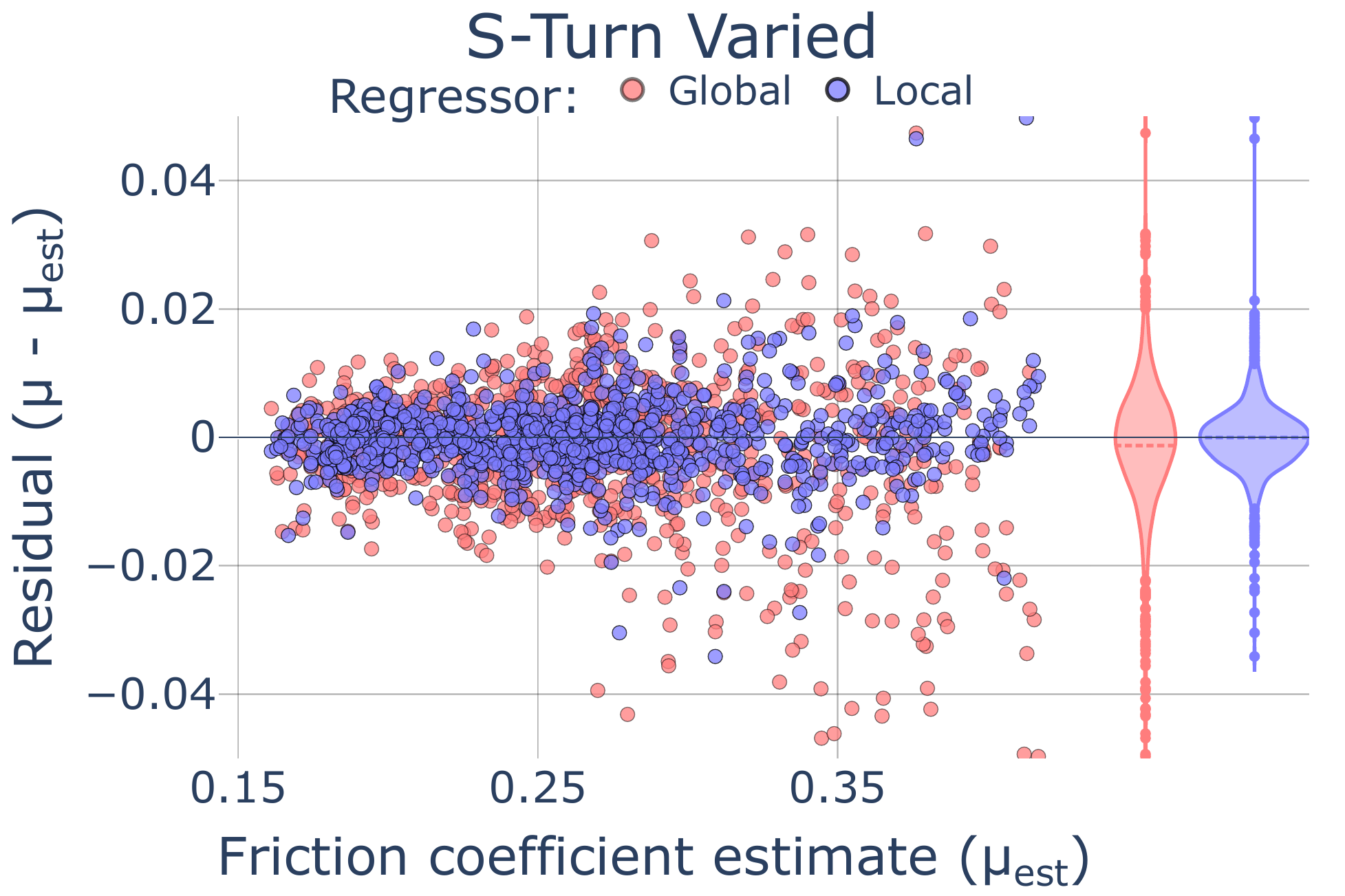}
    \end{subfigure}
\caption{Residual plots of the local and global regressor performance on the test data of the more challenging scenarios}
\label{fig:residualPlots}
\end{figure*}

\begin{figure*}
\includegraphics[width=0.85\textwidth, height=.28\textheight, keepaspectratio]{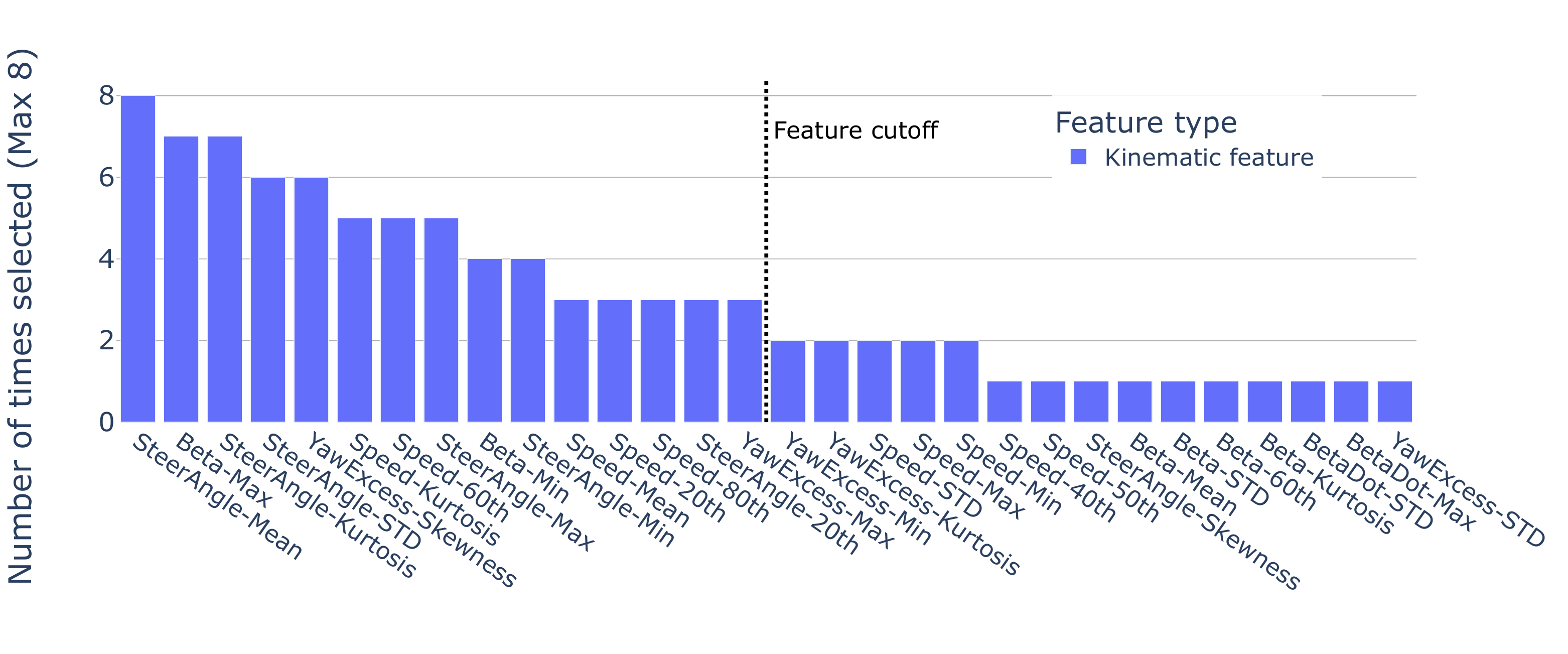}
\includegraphics[width=0.85\textwidth, height=.28\textheight, keepaspectratio]{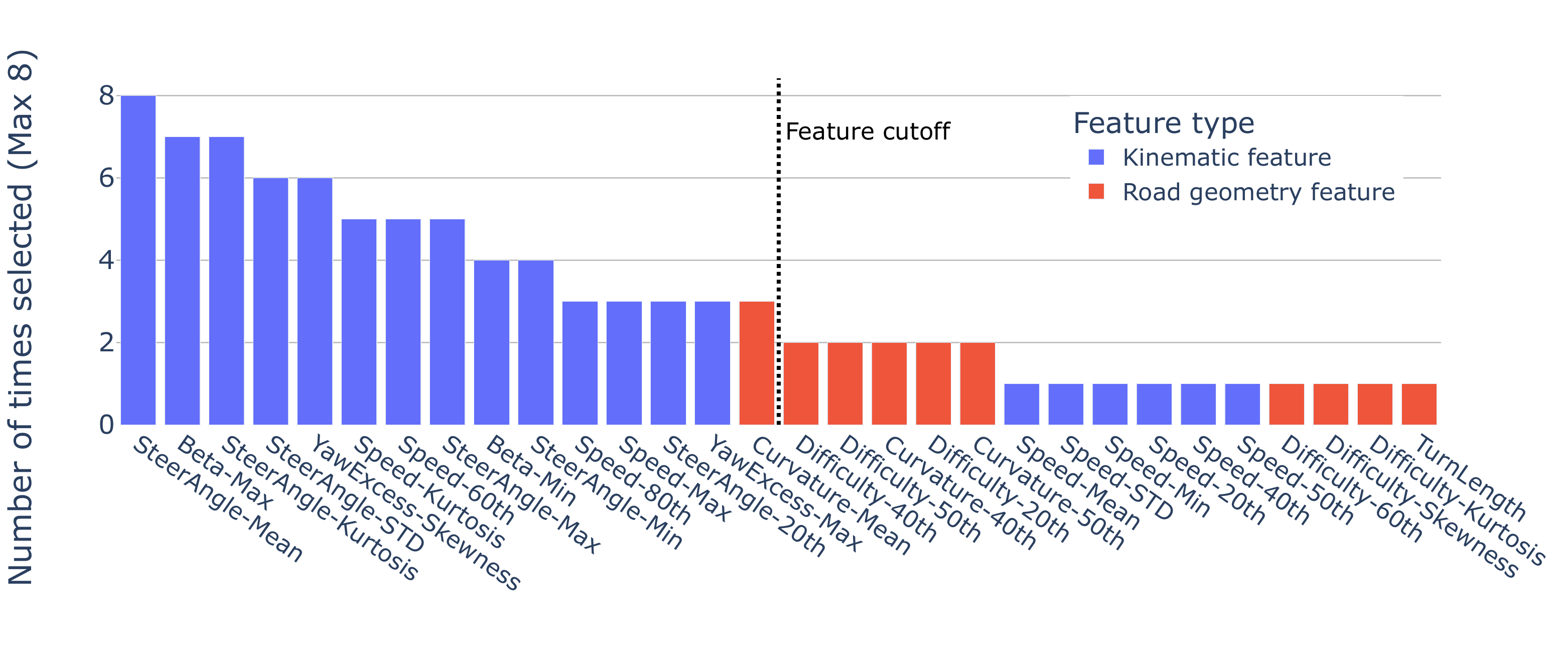}
\caption{The result of running sequential feature selection on input variations - excluding and including road geometry features. 
The features to the left of the vertical dotted line are the ones selected to be used by the regressors.}
\label{fig:selectedFeatures}
\end{figure*}





\end{document}